\documentclass[a4paper,12pt,twoside]{article}
\usepackage{amsmath,amssymb,ifthen}
\usepackage[amsmath,thmmarks]{ntheorem}
\usepackage{paper}
\usepackage[all]{xy}
\usepackage{mathrsfs}
\usepackage{url}
\DeclareMathOperator{\Spf}{Spf}
\DeclareMathOperator{\Spa}{Spa}
\DeclareMathOperator{\rec}{rec}
\DeclareMathOperator{\Ind}{Ind}
\DeclareMathOperator{\Nr}{Nr}
\DeclareMathOperator{\Nrd}{Nrd}
\DeclareMathOperator{\JL}{JL}
\DeclareMathOperator{\SL}{SL}
\DeclareMathOperator{\Iw}{Iw}
\DeclareMathOperator{\cInd}{c-Ind}
\DeclareMathOperator{\Int}{Int}
\DeclareMathOperator{\Lie}{Lie}
\DeclareMathOperator{\Art}{Art}
\DeclareMathOperator{\sgn}{sgn}

\SelectTips{cm}{11}
\title{Parity of the Langlands parameters of conjugate self-dual representations of $\mathrm{GL}(n)$
and the local Jacquet-Langlands correspondence}
\author{Yoichi Mieda}

\def\headertitle{Parity of the Langlands parameters of conjugate self-dual representations}

\begin{document}
\maketitle

\begin{firstfootnote}
 Graduate School of Mathematical Sciences, The University of Tokyo, 3--8--1 Komaba, Meguro-ku, Tokyo, 153--8914, Japan

 E-mail address: \texttt{mieda@ms.u-tokyo.ac.jp}

 2010 \textit{Mathematics Subject Classification}.
 Primary: 11F70;
 Secondary: 11G25, 22E50.
\end{firstfootnote}

\begin{abstract}
 We determine the parity of the Langlands parameter of a conjugate self-dual supercuspidal representation
 of $\GL(n)$ over a non-archimedean local field by means of the local Jacquet-Langlands correspondence. 
 It gives a partial generalization of a previous result on the self-dual case by Prasad and Ramakrishnan.
\end{abstract}

\section{Introduction}

Let $F$ be a $p$-adic field. By the local Langlands correspondence, irreducible smooth representations of
$\GL_n(F)$ are known to be parameterized by $n$-dimensional representations of $W_F\times \SL_2(\C)$,
where $W_F$ denotes the Weil group of $F$. For an irreducible smooth representation $\pi$ of $\GL_n(F)$,
we write $\rec_F(\pi)$ for the attached parameter, which is called the Langlands parameter of $\pi$.

Let us assume that $\pi$ is self-dual, namely, $\pi$ is isomorphic to its contragredient $\pi^\vee$.
Since $\rec_F$ is compatible with dual, $\rec_F(\pi)$ is again self-dual.
Therefore we can consider the problem whether $\rec_F(\pi)$ is symplectic or orthogonal, under the condition
that $\rec_F(\pi)$ is irreducible, in other words, $\pi$ is a discrete series representation.
In \cite{MR2931222}, Prasad and Ramakrishnan answered this question by means of the local Jacquet-Langlands
correspondence.
Let $D$ be a central division algebra of rank $n$ over $F$. Recall that the local Jacquet-Langlands correspondence
(\cite{MR700135}, \cite{MR771672})
gives a bijection between isomorphism classes of irreducible discrete series representations of $\GL_n(F)$
and those of irreducible smooth representations of $D^\times$.
We write $\JL(\pi)$ for the representation of $D^\times$ attached to $\pi$
by this correspondence. The theorem of Prasad and Ramakrishnan is as follows:

\begin{thm}[{{\cite[Theorem B]{MR2931222}}}]\label{thm:P-R}
 Assume that $\pi$ is self-dual. If $n$ is odd, $\rec_F(\pi)$ is always orthogonal
 (this part is clear). If $n$ is even, then $\rec_F(\pi)$ is symplectic (resp.\ orthogonal)
 if and only if $\JL(\pi)$ is orthogonal (resp.\ symplectic).
\end{thm} 

The purpose of this paper is to extend this theorem to the conjugate self-dual setting.
Let $F/F^+$ be a quadratic extension of $p$-adic fields and $\tau$ the generator of $\Gal(F/F^+)$.
A smooth representation $(\pi,V)$ of $\GL_n(F)$ is said to be conjugate self-dual if $\pi^\tau\cong \pi^\vee$,
where $\pi^\tau$ denotes the representation $\GL_n(F)\xrightarrow{\tau}\GL_n(F)\xrightarrow{\pi}\GL(V)$.
If $\pi$ is conjugate self-dual, its Langlands parameter $\rec_F(\pi)$ is also conjugate self-dual
in the following sense. Take $c\in W_{F^+}\setminus W_F$. For a representation $\phi$ of $W_F\times\SL_2(\C)$,
define a new representation $\phi^c$ by $\phi^c(w)=\phi(cwc^{-1})$; it is independent of the choice of $c$
up to isomorphism. A representation $\phi$ is said to be conjugate self-dual if $\phi^c\cong \phi^\vee$ holds.
For an irreducible conjugate self-dual representation $\phi$ of $W_F\times\SL_2(\C)$,
we can define its parity $C_\phi\in\{\pm 1\}$ in the similar way as in the self-dual case
(for the detail, see \cite[\S 3]{MR3202556}, \cite[\S 2.2]{MR3338302} and Section \ref{sec:parity} of this paper).
If $C_\phi=1$, $\phi$ is said to be conjugate orthogonal, otherwise conjugate symplectic.
For an irreducible conjugate self-dual discrete series representation $\pi$,
the parity of $\rec_F(\pi)$ knows whether $\pi$ comes from the standard base change lifting 
or the twisted base change lifting from the quasi-split unitary group $U_{F/F^+}(n)$
(see \cite[\S 2]{MR3338302}).

In this paper, we determine the parity of $\rec_F(\pi)$ by means of $\JL(\pi)$, under the conditions that
\begin{itemize}
 \item $F/F^+$ is at worst tamely ramified,
 \item the invariant of $D$ is $1/n$,
 \item and $\pi$ is supercuspidal (in other words, $\rec_F(\pi)$ is trivial on the $\SL_2(\C)$-factor).
\end{itemize}
Under the first two assumptions, we construct explicitly an automorphism $\tau\colon D^\times\to D^\times$
such that $\tau\vert_{F^\times}$ coincides with $\tau\in\Gal(F/F^+)$, 
and $t\in D^\times$ such that $\tau^2(d)=tdt^{-1}$ for $d\in D^\times$ (Definition \ref{defn:tau-on-D}).
For such a pair $(\tau,t)$, we can define the conjugate self-duality and the parity of 
an irreducible smooth representation of $D^\times$ (see Section \ref{sec:parity}).
Our main theorem is summarized as follows: 

\begin{thm}[Main theorem, Theorem \ref{thm:main}]\label{thm:main-intro}
 Assume that $F/F^+$ is at worst tamely ramified and the invariant of $D$ is $1/n$.
 Let $\pi$ be an irreducible conjugate self-dual supercuspidal representation of $\GL_n(F)$.
 Then, $\JL(\pi)$ is conjugate self-dual with respect to $(\tau,t)$, and its parity $C_{\JL(\pi)}$ satisfies
 \[
  C_{\rec_F(\pi)}=(-1)^{n-1}C_{\JL(\pi)}.
 \]
\end{thm}

Theorems \ref{thm:P-R}, \ref{thm:main-intro} are useful in the study of $\rec_F(\pi)$,
because the determination of $\JL(\pi)$ is usually much easier than that of $\rec_F(\pi)$.
In Section \ref{sec:simple-supercuspidal}, we apply Theorems \ref{thm:P-R}, \ref{thm:main-intro}
to compute the parity of $\rec_F(\pi)$ for conjugate (or usual) self-dual simple supercuspidal representations
of $\GL_n(F)$ (for simple supercuspidal representations,
see \cite{MR2730575}, \cite{MR3164986}, \cite{Knightly-Li}).
For example, we prove that the Langlands parameter of a self-dual simple supercuspidal representation
of $\GL_{2n}(F)$ is symplectic if and only if its central character is trivial.
This result plays a crucial role in the recent study of Oi \cite{Oi-ssc-SO} on the endoscopic lifting of
simple supercuspidal representations of $\mathrm{SO}_{2n+1}(F)$ to $\GL_{2n}(F)$.

Let us explain the strategy of our proof of Theorem \ref{thm:main-intro}.
We use a geometric method.
The non-abelian Lubin-Tate theory (\cite{MR1044827}, \cite{MR1719811}, \cite{MR1876802})
tells us that the correspondences $\rec_F$ and $\JL$
for supercuspidal representations appear in the $\ell$-adic \'etale cohomology of the Lubin-Tate tower,
which is a projective system of universal deformation spaces of
a one-dimensional formal $\mathcal{O}_F$-module $\mathbb{X}$ of height $n$ with suitable level structures.
By using the cup product of the cohomology and a result in \cite{non-cusp}, we can construct
a perfect pairing 
\[
 (\JL(\pi)\boxtimes \rec_F(\pi))\times (\JL(\pi^\vee)\boxtimes \rec_F(\pi^\vee))\to \C
\]
for an irreducible supercuspidal representation $\pi$ of $\GL_n(F)$.
It enables us to compare the parity of $\rec_F(\pi)$ and that of $\JL(\pi)$, provided that $\pi$ is self-dual.
As in the introduction of \cite{MR2931222}, this method had already been found by Fargues;
he announced the supercuspidal case of Theorem \ref{thm:P-R} in \cite[\S 5]{Fargues-Zel} without proof.
The new point of this paper is to adapt the argument above to the conjugate self-dual case.
In the conjugate self-dual case, we need to make the pairing ``Hermitian''.
For this purpose, we introduce a new operator on the Lubin-Tate tower, which we call the twisting operator.
In the definition of it, we need to fix an additional structure on the fixed formal $\mathcal{O}_F$-module
$\mathbb{X}$. This extra structure naturally induces the pair $(\tau,t)$ in Theorem \ref{thm:main-intro},
as $D^\times$ can be identified with the group of self-$\mathcal{O}_F$-isogenies of $\mathbb{X}$.

Since our method is geometric, our theorem is also valid in the equal characteristic case.
On the other hand, we need to assume that the invariant of $D$ is $1/n$ and $\pi$ is supercuspidal,
because this is the only case in which $\rec_F(\pi)$ and $\JL(\pi)$ have nice geometric descriptions.
The author expects that Theorem \ref{thm:main-intro} is true for
any conjugate self-dual discrete series representation $\pi$;
in fact, we can easily verify it for a character twist
of the Steinberg representation (see Remark \ref{rem:Steinberg}).
It seems also an interesting question to extend Theorem \ref{thm:main-intro} to general division algebras.
These problems will be considered in our future works. 

The outline of this paper is as follows. 
In Section \ref{sec:parity}, we give some basic definitions on conjugate self-dual representations
and their parity. We need a slightly general framework than usual,
in order to formulate Theorem \ref{thm:main-intro}. 
Section \ref{sec:proof-main-thm} is devoted to a proof of the main theorem.
After a brief review of the non-abelian Lubin-Tate theory, we introduce and study the twisting operator,
which is a key of our proof. To describe the pair $(\tau,t)$ explicitly, we also need some explicit computations
of Dieudonn\'e modules.
In Section \ref{sec:simple-supercuspidal}, we apply the main theorem to determine the parity of
conjugate self-dual simple supercuspidal representations of $\GL_n(F)$.

\medbreak
\noindent{\bfseries Acknowledgment}\quad
This work was supported by JSPS KAKENHI Grant Number 15K13424.

\medbreak
\noindent{\bfseries Notation}\quad
For a field $L$ and an integer $m\ge 1$, we write $\mu_m(L)$ for the set of $m$th roots of unity in $L$.
If $L$ is a discrete valuation field, we denote the ring of integers of $L$ by $\mathcal{O}_L$,
and the maximal ideal of $\mathcal{O}_L$ by $\mathfrak{p}_L$.
Every representation is considered over $\C$.

\section{Parity of conjugate self-dual representations}\label{sec:parity}
\subsection{Basic definitions and properties}\label{subsec:basic}
Let $G$ be a totally disconnected locally compact topological group.
We fix a continuous automorphism $\tau\colon G\to G$ and an element $t\in G$
satisfying
\[
 \tau^2=\Int(t),\quad \tau(t)=t,
\]
where $\Int(t)\colon G\to G$ denotes the isomorphism $g\mapsto tgt^{-1}$.
For a smooth representation $(\pi,V)$ of $G$, we write $(\pi^\tau,V)$ for the smooth representation
defined by $\pi^\tau(g)=\pi(\tau(g))$.
We say that $\pi$ is conjugate self-dual with respect to $\tau$ if $\pi^\tau$ is isomorphic to
the contragredient representation $\pi^\vee$. If $\pi$ is conjugate self-dual with respect to $\tau$,
we have $\pi^{\vee\vee}\cong (\pi^\tau)^{\vee}=(\pi^\vee)^{\tau}\cong (\pi^\tau)^\tau=\pi^t\cong \pi$
(the last isomorphism is given by $\pi(t)^{-1}$). Hence $\pi$ is admissible.

Let $\pi$ be a smooth representation of $G$ which is conjugate self-dual with respect to $\tau$.
Then, there exists a non-degenerate bilinear pairing
$\langle\ ,\ \rangle\colon V\times V\to \C$ satisfying 
$\langle \pi(\tau(g))x,\pi(g)y\rangle=\langle x,y\rangle$ for every $g\in G$ and $x,y\in V$.
If $\pi$ is irreducible, such a pairing is unique up to scalar by Schur's lemma
(recall that $\pi$ is admissible).

\begin{lem}\label{lem:parity}
 There exists $C_\pi\in\{\pm 1\}$ such that
 $\langle \pi(t)y,x\rangle=C_\pi\langle x,y\rangle$ for every $x,y\in V$.
\end{lem}

\begin{prf}
 Put $\langle x,y\rangle'=\langle \pi(t)y,x\rangle$. 
 Let $g\in G$ and $x,y\in V$ be arbitrary elements, and we put $g'=\tau^{-1}(g)$.
 Then we have
 \begin{align*}
  \langle \pi(\tau(g))x,\pi(g)y\rangle'
  &=\langle \pi(t)\pi(g)y,\pi(\tau(g))x\rangle=\langle \pi(\tau(tg'))y,\pi(tg't^{-1})x\rangle\\
  &=\langle y,\pi(t)^{-1}x\rangle=\langle \pi(t)y,x\rangle=\langle x,y\rangle'.
 \end{align*}
 Therefore there exists $C_\pi\in \C^\times$ such that 
 $\langle x,y\rangle'=C_\pi\langle x,y\rangle$ for every $x,y\in V$.
 
 For $x,y\in V$, we have
 \[
  \langle x,y\rangle=\langle\pi(\tau(t))x,\pi(t)y\rangle=\langle \pi(t)x,\pi(t)y\rangle
 =C_\pi\langle \pi(t)y,x\rangle=C_\pi^2\langle x,y\rangle.
 \]
 Hence we have $C_\pi^2=1$. This concludes the proof.
\end{prf}

\begin{rem}
 The sign $C_\pi$ depends not only on $\tau$ but also on $t$. Let $t'\in G$ be another element
 satisfying $\tau^2=\Int(t')$. Then $z=t't^{-1}$ lies in the center of $G$
 and fixed by $\tau$.
 It is immediate to see that $C_\pi$ for $t'$ equals $\omega_\pi(z)C_\pi$, where $\omega_\pi$ denotes
 the central character of $\pi$.
\end{rem}

We call $C_\pi$ the parity of $\pi$ (with respect to $(\tau,t)$).
If $C_\pi=1$ (resp.\ $C_\pi=-1$), we say that $\pi$ is conjugate orthogonal
(resp.\ conjugate symplectic). If $\tau=\id$ and $t=1$, this notion coincides with the standard one.

\begin{rem}\label{rem:parity-odd-dim}
 Consider the case where $(\pi,V)$ is finite-dimensional, and put $m=\dim_\C V$.
 \begin{enumerate}
  \item Assume that $m=1$, and identify $V$ with $\C$. 
	Then, $\langle\ ,\ \rangle\colon \C\times \C\to \C$; $(x,y)\mapsto xy$ gives
	a non-degenerate bilinear pairing satisfying $\langle \pi(\tau(g))x,\pi(g)y\rangle=\langle x,y\rangle$.
	From this pairing we can deduce $C_\pi=\pi(t)$.
  \item Let $\langle\ ,\ \rangle\colon V\times V\to \C$ be a non-degenerate bilinear pairing
	as in the definition of the parity. Put $(\det\pi,\det V)=(\bigwedge^m \pi,\bigwedge^m V)$.
	Then, $\langle\ ,\ \rangle$ induces a pairing $\det V\times\det V\to \C$ by
	\[
	 \langle x_1\wedge\cdots\wedge x_m,y_1\wedge\cdots\wedge y_m\rangle=\sum_{\sigma\in\mathfrak{S}_m}\sgn(\sigma)\langle x_1,y_{\sigma(1)}\rangle\cdots \langle x_m,y_{\sigma(m)}\rangle.
	\]
	It is non-degenerate and satisfies
	\[
	 \langle (\det\pi)(\tau(g))x,(\det\pi)(g)y\rangle=\langle x,y\rangle, \quad \langle (\det\pi)(t)y,x\rangle=C_\pi^m\langle x,y\rangle
	\]
	for $x,y\in\det V$ and $g\in G$. Hence we have $C_{\det\pi}=C_\pi^m$.
 \end{enumerate}
 In particular, if $m$ is odd, the parity $C_\pi$ can be computed as follows:
 \[
  C_\pi=C_\pi^m=C_{\det\pi}=\det \pi(t).
 \]
 In contrast, if $m$ is even, the parity is a more subtle invariant.
\end{rem}

We give two elementary lemmas.

\begin{lem}\label{lem:product}
 Assume that $(G,\tau,t)$ is decomposed into $(G_1\times G_2,\tau_1\times \tau_2,(t_1,t_2))$,
 where $G_i$ is a totally disconnected locally compact topological group,
 $\tau_i\colon G_i\to G_i$ a continuous automorphism and $t_i\in G_i$ satisfying $\tau_i^2=\Int(t_i)$.
 For each $i=1,2$, let $(\pi_i,V_i)$ be an irreducible smooth representation of $G_i$
 conjugate self-dual with respect to $\tau_i$.
 Then, $(\pi_1\boxtimes \pi_2,V_1\otimes V_2)$ is an irreducible smooth representation of $G$
 conjugate self-dual with respect to $\tau$, and $C_{\pi_1\boxtimes \pi_2}$ is equal to $C_{\pi_1}C_{\pi_2}$.
\end{lem}

\begin{prf}
 It is well-known that the exterior tensor product of irreducible admissible representations
 is irreducible. The parity can be computed by using the pairing
 $\langle x_1\otimes x_2,y_1\otimes y_2\rangle=\langle x_1,y_1\rangle_1\langle x_2,y_2\rangle_2$,
 where $\langle\ ,\ \rangle_i\colon V_i\times V_i\to \C$ is an appropriate pairing attached to $\pi_i$.
\end{prf}

\begin{lem}\label{lem:change-tau}
 Take an element $h\in G$ and put $\tau'=\Int(h)\circ\tau$, $t'=h\tau(h)t$.
 Then we have $\tau'^2=\Int(t')$. For an irreducible smooth representation $\pi$ of $G$,
 $\pi$ is conjugate self-dual with respect to $\tau$ if and only if
 it is conjugate self-dual with respect to $\tau'$.
 If $\pi$ is conjugate self-dual with respect to $\tau$ and $\tau'$, its parity with respect to $(\tau,t)$
 coincides with that with respect to $(\tau',t')$.
\end{lem}

\begin{prf}
 The claim $\tau'^2=\Int(t')$ is immediate. We write $V$ for the representation space of $\pi$.
 Assume that $\pi$ is conjugate self-dual with respect to $\tau$, and take a non-degenerate pairing
 $\langle\ ,\ \rangle\colon V\times V\to \C$ satisfying $\langle \pi(\tau(g))x,\pi(g)y\rangle=\langle x,y\rangle$.
 Let $\langle\ ,\ \rangle_h\colon V\times V\to \C$ be the pairing defined by 
 $\langle x,y\rangle_h=\langle \pi(h)^{-1}x,y\rangle$.
 It is a non-degenerate pairing and satisfies
 \begin{align*}
  \langle \pi(\tau'(g))x,\pi(g)y\rangle_h&=\langle \pi(h)^{-1}\pi(h\tau(g)h^{-1})x,\pi(g)y\rangle
  =\langle \pi(\tau(g))\pi(h)^{-1}x,\pi(g)y\rangle\\
  &=\langle \pi(h)^{-1}x,y\rangle=\langle x,y\rangle_h.
 \end{align*}
 Therefore $\pi^{\tau'}\cong\pi^\vee$, that is, $\pi$ is conjugate self-dual with respect to $\tau'$.
 Since $\tau=\Int(h^{-1})\circ\tau'$, the converse is also the case.

 Let us denote by $C$ (resp.\ $C'$) the parity of $\pi$, which is assumed to be conjugate self-dual,
 with respect to $(\tau,t)$ (resp.\ $(\tau',t')$). We use the pairing $\langle\ ,\ \rangle_h$ to compute $C'$.
 For $x,y\in V$, we have
 \begin{align*}
 C'\langle x,y\rangle_h&=\langle \pi(t')y,x\rangle_h
  =\langle \pi(h^{-1}t')y,x\rangle=\langle \pi(\tau(h)t)y,x\rangle
  =\langle \pi(t)y,\pi(h)^{-1}x\rangle\\
  &=C\langle \pi(h)^{-1}x,y\rangle=C\langle x,y\rangle_h.
 \end{align*}
 Hence we conclude that $C=C'$.
\end{prf}

Let $H$ be an open subgroup of $G$ such that $H\backslash G$ is compact.
Take a smooth character $\chi\colon H\to \C^\times$ such that $\pi=\cInd_H^G\chi$ is irreducible.
We will consider when $\pi$ is conjugate self-dual with respect to $\tau$, and
how to compute the parity of $\pi$.

\begin{prop}\label{prop:parity-cInd}
 Put $H^\tau=\tau^{-1}(H)$, and write $\chi^\tau$ for the character
 $H^\tau\to \C^\times$; $h\mapsto \chi(\tau(h))$.
 Assume that there exists $a\in G$ which intertwines $(H,\chi^{-1})$ and $(H^\tau,\chi^\tau)$; namely,
 satisfies the following conditions:
 \[
  aHa^{-1}=H^\tau,\quad \chi(h)^{-1} =\chi^\tau(aha^{-1})\ \text{for every $h\in H$}.
 \]
 Then, the representation $\pi=\cInd_H^G\chi$ is conjugate self-dual with respect to $\tau$.
 Furthermore, an element $z=\tau(a)ta$ lies in $H$, and the parity $C_\pi$ of $\pi$ is given by $\chi(z)$.
\end{prop} 

\begin{prf}
 For $f\in \cInd_H^G\chi$, let $f^\tau\colon G\to \C$ be the function $g\mapsto f(\tau(g))$.
 Then, it is easy to see that $f^\tau$ belongs to $\cInd_{H^\tau}^G\chi^\tau$, and $f\mapsto f^\tau$
 gives an isomorphism $(\cInd_H^G\chi)^\tau\xrightarrow{\cong}\cInd_{H^\tau}^G\chi^\tau$ of $G$-representations.
 On the other hand, for $f\in \cInd_{H^\tau}^G\chi^{\tau}$, 
 let $f^a\colon G\to \C$ be the function $g\mapsto f(ag)$.
 We can check that $f^a$ belongs to $\cInd_H^G\chi^{-1}$ and $f\mapsto f^a$ gives an isomorphism
 $\cInd_{H^\tau}^G\chi^{\tau}\xrightarrow{\cong}\cInd_H^G\chi^{-1}$.
 Since $H$ is cocompact in $G$, we have $(\cInd_H^G\chi)^\vee\cong \Ind_H^G\chi^{-1}=\cInd_H^G\chi^{-1}$.
 Hence we have $\pi^\tau=(\cInd_H^G\chi)^\tau\cong \cInd_{H^\tau}^G\chi^\tau\cong \cInd_H^G\chi^{-1}\cong \pi^\vee$.
 In other words, $\pi$ is conjugate self-dual with respect to $\tau$.

 Next we prove $z\in H$. First we will see that $z$ normalizes $(H,\chi)$.
 Since $H^\tau=aHa^{-1}$, we have $H=\tau(a)\tau^2(H^\tau)\tau(a)^{-1}=\tau(a)tH^\tau t^{-1}\tau(a)^{-1}=zHz^{-1}$.
 Therefore $z$ normalizes $H$. Moreover, for $h\in H$ we have
 \begin{align*}
  \chi(z^{-1}hz)&=\chi^\tau(az^{-1}hza^{-1})^{-1}=\chi^\tau(t^{-1}\tau(a)^{-1}h\tau(a)t)^{-1}
  =\chi(a^{-1}t^{-1}\tau(h)ta)^{-1}\\
  &=\chi^\tau(t^{-1}\tau(h)t)=\chi(h).
 \end{align*}
 Thus $z$ fixes $\chi$.

 Recall that we are assuming that $\pi=\cInd_H^G\chi$ is irreducible. Therefore,
 \[
  \Hom_G(\pi,\pi)=\Hom_H(\chi,(\cInd_H^G\chi)\vert_H)\cong \Hom_H\Bigl(\chi,\bigoplus_{g\in H\backslash G/H}\cInd_{H\cap g^{-1}Hg}^H\chi^g\Bigr)
 \]
 is one-dimensional (here $\chi^g$ denotes the character $h'\mapsto \chi(gh'g^{-1})$ on $H\cap g^{-1}Hg$).
 Since $\cInd_{H\cap z^{-1}Hz}^H\chi^z=\chi$, $z$ must lie in $H$; otherwise the direct sum above contains
 $\chi\oplus \chi$. 

 Finally we compute the parity of $\pi$.
 Define a pairing $\langle\ ,\ \rangle\colon \cInd_H^G\chi\times \cInd_H^G\chi\to \C$ by
 \[
  \langle f_1,f_2\rangle=\sum_{g\in H\backslash G}(f_1^\tau)^a(g)f_2(g).
 \]
 This is the composite of
 $(\cInd_H^G\chi)^\tau\times \cInd_H^G\chi\xrightarrow[\cong]{(1)}\cInd_H^G\chi^{-1}\times \cInd_H^G\chi\xrightarrow{(2)}\C$, where (1) denotes the isomorphism $(f_1,f_2)\mapsto ((f_1^\tau)^a,f_2)$ and (2) the canonical pairing.
 Hence $\langle\ ,\ \rangle$ is a non-degenerate pairing satisfying 
 $\langle \pi(\tau(g))f_1,\pi(g)f_2\rangle=\langle f_1,f_2\rangle$ for every $g\in G$ and
 $f_1,f_2\in \cInd_H^G\chi$.

 By definition we can compute as follows:
 \begin{align*}
  \langle \pi(t)f_2,f_1\rangle&=\sum_{g\in H\backslash G}(f_2^\tau)^a(gt)f_1(g)
 =\sum_{g\in H\backslash G}f_2(\tau(agt))f_1(g)\\
  &\stackrel{(*)}{=}\sum_{g'\in H\backslash G}f_2(g')f_1(a^{-1}t^{-1}\tau(g'))
  =\sum_{g'\in H\backslash G}f_1(z^{-1}\tau(ag'))f_2(g')\\
  &=\sum_{g'\in H\backslash G}\chi(z)^{-1}f_1(\tau(ag'))f_2(g')=\chi(z)^{-1}\langle f_1,f_2\rangle.
 \end{align*}
 At the equality $(*)$, we put $g'=\tau(agt)$. As $\tau(aHgt)=\tau(aHa^{-1})\tau(agt)=H\tau(agt)$,
 this replacement is well-defined.
 Hence the parity $C_\pi=C_\pi^{-1}$ of $\pi$ equals $\chi(z)$. This completes the proof.
\end{prf}

Proposition \ref{prop:parity-cInd} will be used in Section \ref{sec:simple-supercuspidal}.

\subsection{Division algebra setting}\label{subsec:div-alg}
Let $F^+$ be a non-archimedean local field and $F$ a separable extension of $F^+$ such that $[F:F^+]\le 2$.
Denote by $\tau$ the generator of $\Gal(F/F^+)$.
Let $q$ (resp.\ $q'$) denote the cardinality of the residue field of $\mathcal{O}_F$ (resp.\ $\mathcal{O}_{F^+}$).
We denote the characteristic of $\F_q$ by $p$.

The extension $F/F^+$ provides two well-known examples of $(G,\tau,t)$ in the previous subsection.

\begin{exa}\label{exa:GL_n}
 For an integer $n\ge 1$, put $G=\GL_n(F)$. Let $\tau\colon G\to G$ be an automorphism induced by
 $\tau\in \Gal(F/F^+)$. Then we have $\tau^2=\id$, and we can set $t=1$.
\end{exa}

\begin{exa}\label{exa:Weil-group}
 Let $G$ be the Weil group $W_F$ of $F$.
 Fix an element $c\in W_{F^+}$ whose image in $W_{F^+}/W_F$ is a generator, and let
 $\tau\colon G\to G$ be $\Int(c)$. Then $c^2$ lies in $W_F$, and we can set $t=c^2$.
 The conjugate self-duality and the parity are independent of the choice of $c$.
 Indeed, another choice of $c$ is of the form $wc$ with $w\in W_F$. 
 Use Lemma \ref{lem:change-tau} to $\tau'=\Int(wc)=\Int(w)\circ \tau$ and $t'=(wc)^2=w(cwc^{-1})c^2=w\tau(w)t$.

 The conjugate self-duality and the parity in this case coincides with those in \cite[\S 3]{MR3202556}
 and \cite[\S 2.2]{MR3338302}.
\end{exa}

The parity under the setting in Example \ref{exa:Weil-group} is interesting because of the following theorem:

\begin{thm}
 Let $\pi$ be an irreducible supercuspidal representation of $\GL_n(F)$ and $\rec_F(\pi)$
 the corresponding $n$-dimensional irreducible smooth representation of $W_F$
 under the local Langlands correspondence.
 \begin{enumerate}
  \item The representation $\pi$ is conjugate self-dual under the setting in Example \ref{exa:GL_n}
	if and only if $\rec_F(\pi)$ is conjugate self-dual under the setting in Example \ref{exa:Weil-group}.
  \item Assume that $F\neq F^+$ and the characteristic of $F$ is $0$.
	The representation $\pi$ belongs to the image of the standard (resp.\ twisted) base change lift
	from the quasi-split unitary group
	$U_{F/F^+}(n)$ if and only if the parity $C_{\rec_F(\pi)}$ is equal to $(-1)^{n-1}$ (resp.\ $(-1)^n$).
	For the notion of the base change lift, see \cite[\S 2]{MR3338302}.
 \end{enumerate}
\end{thm}

In the following, we introduce another setting. Fix a separable closure $\overline{F}$ and 
a uniformizer $\varpi$ of $F$.
For an integer $n\ge 1$, we denote by $F_n$ (resp.\ $F_n^+$) the unramified extension of degree $n$ of $F$
(resp.\ $F^+$) contained in $\overline{F}$, and by $\sigma\in \Gal(F_n/F)$ the arithmetic Frobenius lift.
Let $D$ be the central division algebra over $F$ with invariant $1/n$.
Recall that $D$ can be written as $F_n[\Pi]$, where $\Pi$ is an element satisfying $\Pi^n=\varpi$ and
$\Pi a=\sigma(a)\Pi$ for every $a\in F_n$.
Assuming the tameness of $F/F^+$, we will explicitly construct an isomorphism $\tau\colon D\to D$
whose restriction to the center $F$ coincides with $\tau\in \Gal(F/F^+)$. 

\begin{defn}\label{defn:tau-on-D}
 Assume that $F/F^+$ is at worst tamely ramified.
 \begin{enumerate}
  \item Suppose that $F/F^+$ is an unramified quadratic extension.
	Then, $\tau\in \Gal(F/F^+)$ is canonically extended to the arithmetic Frobenius lift in $\Gal(F_n/F^+)$,
	that is also denoted by $\tau$. It satisfies $\sigma=\tau^2$.
	In this case, we take $\varpi$ in $F^+$ and define $\tau\colon D\to D$ by
	$a\mapsto \tau(a)$ $(a\in F_n)$ and $\Pi\mapsto \Pi$. We put $t=\Pi$.
  \item Suppose that $F/F^+$ is a ramified quadratic extension (thus $p\neq 2$).
	Then, the restriction map $\Gal(F_n/F_n^+)\to \Gal(F/F^+)$ is an isomorphism.
	We also write $\tau$ for the generator of $\Gal(F_n/F_n^+)$.
	It commutes with $\sigma\in \Gal(F_n/F)$.

	In this case, we can (and do) take $\varpi$ so that $\tau(\varpi)=-\varpi$.
	Fix an element $\beta\in \overline{F}$ such that $\beta^{q^n-1}=-1$
	and put $\alpha=\beta^{q-1}$.
	Since $\alpha^{q^n-1}=(-1)^{q-1}=1$, $\alpha$ belongs to 
	$\mu_{q^n-1}(\overline{F})=\mu_{q^n-1}(F_n^+)$ and
	$\Nr_{F_n/F}(\alpha)=\alpha^{1+q+\cdots+q^{n-1}}=\beta^{q^n-1}=-1$.
	We define $\tau\colon D\to D$ by $a\mapsto \tau(a)$ $(a\in F_n)$ and $\Pi\mapsto \alpha\Pi$.
	Note that $(\alpha\Pi)^n=\Nr_{F_n/F}(\alpha)\Pi^n=-\varpi=\tau(\varpi)$
	and $(\alpha\Pi)\tau(a)=\alpha\sigma(\tau(a))\Pi=\tau(\sigma(a))(\alpha\Pi)$, 
	which ensure the well-definedness of $\tau$.
	We put $t=\beta^{-2}\in \mu_{q^n-1}(\overline{F})=\mu_{q^n-1}(F_n^+)$.
  \item If $F=F^+$, then we define $\tau\colon D\to D$ to be the identity map.
	We put $t=1$. 
 \end{enumerate}
\end{defn}

In each case we can check that $\tau^2(d)=tdt^{-1}$ holds for every $d\in D$.
Therefore, the triple $(D^\times,\tau,t)$ gives an example of the setting in Section \ref{subsec:basic}.

\begin{rem}\label{rem:tau-on-D}
 \begin{enumerate}
  \item In the second case, the conjugate self-duality and the parity are independent of the choice of $\beta$.
	Indeed, let $\beta'\in\overline{F}$ be another element such that $\beta'^{q^n-1}=-1$.
	Then $\gamma=\beta/\beta'$ lies in $\mu_{q^n-1}(\overline{F})=\mu_{q^n-1}(F_n^+)$.
	We put $\alpha'=\beta'^{q-1}$ and write $\tau'$, $t'$ for $\tau$, $t$ attached to $\beta'$, respectively.
	Since $\alpha'\Pi=\gamma(\alpha\Pi)\gamma^{-1}$, we have $\tau'=\Int(\gamma)\circ\tau$
	and $t'=\beta'^{-2}=\gamma^2t=\gamma\tau(\gamma)t$. Hence the independence follows from
	Lemma \ref{lem:change-tau}.
  \item In the second case, assume that $n$ is odd. Take $\varepsilon\in \F_q^\times\setminus (\F_q^\times)^2$
	and $\eta\in \F_{q^2}$ such that $\eta^2=\varepsilon^{-1}$.
	We have $\eta^{q-1}=-1$.
	
	Then, the unique lifting $\beta\in \mu_{q^2-1}(\mathcal{O}_{F_2})$ of $\eta$ satisfies
	$\beta^{q^n-1}=(-1)^{1+q+\cdots+q^{n-1}}=-1$. 
	Under this choice of $\beta$, we have $\alpha=\beta^{q-1}=-1$. Moreover, the element $t=\beta^{-2}$
	is the unique element of $\mu_{q-1}(\mathcal{O}_F)$ lifting $\varepsilon$.
 \end{enumerate}
\end{rem}

Our main theorem is as follows:

\begin{thm}\label{thm:main}
 Assume that $F/F^+$ is at worst tamely ramified.
 Let $\pi$ be an irreducible supercuspidal representation of $\GL_n(F)$
 which is conjugate self-dual under the setting in Example \ref{exa:GL_n}.
 We write $\JL(\pi)$ for the irreducible smooth representation of $D^\times$ attached to $\pi$
 under the local Jacquet-Langlands correspondence.

 Then, $\JL(\pi)$ is conjugate self-dual with respect to $\tau\colon D^\times\to D^\times$
 introduced in Definition \ref{defn:tau-on-D}. Moreover, we have
 \[
  C_{\rec_F(\pi)}=(-1)^{n-1}C_{\JL(\pi)}.
 \]
\end{thm}

\begin{rem}\label{rem:Steinberg}
 \begin{enumerate}
  \item The case where $F=F^+$ and the characteristic of $F$ is $0$ has been obtained in \cite{MR2931222},
	in which a discrete series representation $\pi$ is also treated.
	The same statement for the case $F=F^+$ is also announced in \cite[\S 5]{Fargues-Zel} without proof.
  \item It is natural to expect that Theorem \ref{thm:main} remains true for
	conjugate self-dual discrete series representations of $\GL_n(F)$.
	For example, let us consider a twist of the Steinberg representation $\pi=\mathbf{St}\otimes(\chi\circ\det)$,
	where $\chi\colon F^\times\to \C^\times$ is a smooth character.
	Since $\mathbf{St}^\tau\cong \mathbf{St}\cong \mathbf{St}^\vee$,
	the representation $\pi$ is conjugate self-dual if and only if $\chi^\tau=\chi^{-1}$.
	The Langlands parameter $\rec_F(\pi)$ is given by 
	$(\chi\circ\Art_F^{-1})\boxtimes \mathrm{Sym}^{n-1}\mathbf{Std}\colon W_F\times \SL_2(\C)\to \GL_n(\C)$, 
	where $\Art_F\colon F^\times\xrightarrow{\cong}W_F^{\mathrm{ab}}$ denotes the isomorphism of
	the local class field theory, and $\mathbf{Std}$ the standard representation of $\SL_2(\C)$.
	The parity of $\mathrm{Sym}^{n-1}\mathbf{Std}$ equals $(-1)^{n-1}$.
	By Remark \ref{rem:parity-odd-dim} (i), the parity of $\chi\circ\Art_F^{-1}$ is given by
	$\chi(\Art_F^{-1}(c^2))=\chi(\Art_{F^+}^{-1}(c))$ (recall that the image of $c$ under the transfer map
	$W_{F^+}^{\mathrm{ab}}\to W_F^{\mathrm{ab}}$ is $c^2$).
	Hence we obtain $C_{\rec_F(\pi)}=(-1)^{n-1}\chi(\Art_{F^+}^{-1}(c))$.
	On the other hand, we have $\JL(\pi)=\chi\circ\Nrd$, where $\Nrd$ denotes the reduced norm of $D$.
	Its parity $C_{\chi\circ\Nrd}$ equals $\chi(\Nrd(t))$.
	By definition, both $\Art_{F^+}^{-1}(c)$ and $\Nrd(t)$ lie in $(F^+)^\times\setminus \Nr_{F/F^+}(F^\times)$.
	Since $\chi\vert_{(F^+)^\times}$ factors through $(F^+)^\times/\Nr_{F/F^+}(F^\times)$,
	we conclude that $\chi(\Art_{F^+}^{-1}(c))=\chi(\Nrd(t))$ and $C_{\rec_F(\pi)}=(-1)^{n-1}C_{\JL(\pi)}$.
 \end{enumerate}
\end{rem}

\section{Proof of the main theorem}\label{sec:proof-main-thm}
\subsection{Review of the non-abelian Lubin-Tate theory}\label{subsec:NALT}
To prove Theorem \ref{thm:main}, we use the non-abelian Lubin-Tate theory, which is a geometric
realization of the local Langlands correspondence for $\GL_n$.
Here we recall it briefly. Let $F$ be a non-archimedean local field and $\varpi$ its uniformizer.
Take an integer $n\ge 1$.
We write $F^{\mathrm{ur}}$ for the maximal unramified extension of $F$
inside the fixed separable closure $\overline{F}$, and $\breve{F}$ for the completion of $F^{\mathrm{ur}}$.

Let $\mathbf{Nilp}$ be the category of schemes over $\mathcal{O}_{\breve{F}}$ on which $\varpi$ is
locally nilpotent. For an object $S$ of $\mathbf{Nilp}$, we denote the structure morphism 
$S\to \Spec\mathcal{O}_{\breve{F}}$ by $\phi_S$. Put $\overline{S}=S\otimes_{\mathcal{O}_{\breve{F}}}\mathcal{O}_{\breve{F}}/\mathfrak{p}_{\breve{F}}$.
Recall that a formal $\mathcal{O}_F$-module over $S$ is a formal group $X$ over $S$ endowed with
an $\mathcal{O}_F$-action $\iota\colon \mathcal{O}_F\to \End(X)$ such that the following two actions of
$\mathcal{O}_F$ on the Lie algebra $\Lie(X)$ coincide:
\begin{itemize}
 \item the action induced by $\iota$, and
 \item that induced by the $\mathcal{O}_S$-module structure of $\Lie(X)$ and the structure homomorphism
       $\mathcal{O}_F\to \mathcal{O}_{\breve{F}}\to \mathcal{O}_S$.
\end{itemize}

Fix a one-dimensional formal $\mathcal{O}_F$-module $\mathbb{X}$ of $\mathcal{O}_F$-height $n$ over 
$\overline{\F}_q=\mathcal{O}_{\breve{F}}/\mathfrak{p}_{\breve{F}}$.
Such $\mathbb{X}$ is unique up to isomorphism.
Put $D=\End_{\mathcal{O}_F}(\mathbb{X})\otimes_{\Z}\Q$, which is known to be a central division algebra over $F$
with invariant $1/n$.

Let $\mathcal{M}\colon \mathbf{Nilp}\to \mathbf{Set}$ be the functor that sends
$S$ to the set of isomorphism classes of pairs $(X,\rho)$, where $X$ is a formal $\mathcal{O}_F$-module
over $S$ and $\rho\colon \phi_{\overline{S}}^*\mathbb{X}\to X\times_S\overline{S}$ is
an $\mathcal{O}_F$-quasi-isogeny.
It is known that $\mathcal{M}$ is represented by a formal scheme over $\mathcal{O}_{\breve{F}}$,
which is non-canonically isomorphic to the disjoint union of countable copies of 
$\Spf \mathcal{O}_{\breve{F}}[[T_1,\ldots,T_{n-1}]]$ (see \cite{MR0238854}, \cite{MR0384707}, \cite{MR1393439}).
The group of self-isogenies $\mathbf{QIsog}_{\mathcal{O}_F}(\mathbb{X})=D^\times$ naturally acts on $\mathcal{M}$
on the right; $h\in D^\times$ sends $(X,\rho)$ to $(X,\rho\circ \phi_{\overline{S}}^*h)$.
The formal scheme $\mathcal{M}$ is endowed with another structure, called a Weil descent datum.
It is an isomorphism $\alpha\colon \mathcal{M}\to \mathcal{M}$ that makes the following diagram commute:
\[
 \xymatrix{%
 \mathcal{M}\ar[r]^-{\alpha}\ar[d]& \mathcal{M}\ar[d]\\
 \Spf\mathcal{O}_{\breve{F}}\ar[r]^-{\sigma^*}& \Spf\mathcal{O}_{\breve{F}}\lefteqn{.}
 }
\]
Here $\sigma\colon \mathcal{O}_{\breve{F}}\to \mathcal{O}_{\breve{F}}$ is induced from the unique element 
$\sigma\in \Gal(F^{\mathrm{ur}}/F)$ lifting the arithmetic Frobenius automorphism
$\overline{\sigma}\colon \overline{\F}_q\to \overline{\F}_q$, as in Section \ref{subsec:div-alg}.
In order to describe this isomorphism, it suffices to construct a bijection 
$\alpha\colon \mathcal{M}(S)\to \mathcal{M}(S^\sigma)$ for each $S\in \mathbf{Nilp}$ compatibly,
where $S^\sigma$ denotes the object 
$S\xrightarrow{\phi_S}\Spec\mathcal{O}_{\breve{F}}\xrightarrow{\sigma^*}\Spec\mathcal{O}_{\breve{F}}$
of $\mathbf{Nilp}$.
For $(X,\rho)\in\mathcal{M}(S)$, we define 
$\alpha(X,\rho)=(X,\rho\circ \phi_{\overline{S}}^*\Frob^{-1}_{\mathbb{X}})$,
where $\Frob_{\mathbb{X}}\colon \mathbb{X}\to (\overline{\sigma}^*)^*\mathbb{X}$ denotes the $q$th power Frobenius
morphism, which is an $\mathcal{O}_F$-isogeny of $\mathcal{O}_F$-height $1$.

Next we consider level structures. For $m\ge 0$, let $\mathcal{M}_m\colon \mathbf{Nilp}\to \mathbf{Set}$
be the functor that sends $S$ to the set of isomorphism classes of triples $(X,\rho,\eta)$,
where $(X,\rho)\in \mathcal{M}(S)$ and $\eta$ is a Drinfeld $m$-level structure on $X$
(for its definition, see \cite[\S 4]{MR0384707} and \cite[\S II.2]{MR1876802}).
It is represented by a formal scheme finite and flat over $\mathcal{M}$,
and $\{\mathcal{M}_m\}_{m\ge 0}$ form a projective system called the Lubin-Tate tower.
The action of $D^\times$ and the Weil descent datum on $\mathcal{M}$ naturally extend to $\mathcal{M}_m$,
and they are compatible with the transition morphisms of the tower.
Further, the group $\GL_n(F)$ acts on $\{\mathcal{M}_m\}_{m\ge 0}$ on the right as a pro-object
(see \cite[\S 2.2]{MR2383890} for the definition).
This action is called the Hecke action. The principal congruence subgroup 
$K_m=\Ker(\GL_n(\mathcal{O}_F)\to \GL_n(\mathcal{O}_F/\mathfrak{p}_F^m))$ of $\GL_n(F)$ acts trivially on $\mathcal{M}_m$.

By taking the rigid generic fiber, we obtain a projective system $\{M_m\}_{m\ge 0}$ of rigid spaces,
whose transition maps are finite and \'etale.
Each $M_m$ is an $n-1$-dimensional smooth rigid space over $\breve{F}$.
For a compact open subgroup $K$ of $\GL_n(\mathcal{O}_F)$, we can define the rigid space $M_K$
as the quotient of $M_m$ by $K/K_m$, where $m\ge 0$ is an integer satisfying $K_m\subset K$.
It is independent of the choice of $m$, and $M_{K_m}$ coincides with $M_m$.
These rigid spaces form a projective system $\{M_K\}_{K\subset \GL_n(\mathcal{O}_F)}$
with finite \'etale transition maps. The actions of $D^\times$ and $\GL_n(F)$, and the Weil descent datum
naturally extend to it.

For a discrete cocompact subgroup $\Gamma$ of $F^\times$ (e.g., $\varpi^{d\Z}$ for an integer $d\ge 1$),
we may consider the
quotient towers $\{\mathcal{M}_m/\Gamma\}_m$ and $\{M_K/\Gamma\}_K$,
where $\Gamma$ is regarded as a discrete subgroup of $D^\times$ by $F^\times\subset D^\times$.
It is known that the actions of $\GL_n(F)$ on these towers
are trivial on $\Gamma\subset F^\times\subset \GL_n(F)$ (see \cite[Lemma 5.36]{MR1393439}).

Now we take a prime number $\ell\neq p$ and consider the $\ell$-adic \'etale cohomology
of the Lubin-Tate tower
\[
 H^i_{\mathrm{LT}/\Gamma,c}=\varinjlim_K H^i_c\bigl((M_K/\Gamma)\otimes_{\breve{F}}\widehat{\overline{F}},\overline{\Q}_\ell\bigr),\quad
 H^i_{\mathrm{LT}/\Gamma}=\varinjlim_K H^i\bigl((M_K/\Gamma)\otimes_{\breve{F}}\widehat{\overline{F}},\overline{\Q}_\ell\bigr),
\]
where $\widehat{\overline{F}}$ denotes the completion of $\overline{F}$.
The groups $\GL_n(F)$ and $D^\times$ act on $H^i_{\mathrm{LT}/\Gamma,c}$ and $H^i_{\mathrm{LT}/\Gamma}$.
The actions of $\GL_n(F)$ on both spaces are obviously smooth, and moreover admissible. 
The action of $D^\times$ on $H^i_{\mathrm{LT}/\Gamma,c}$ is also known to be smooth
(see \cite[Lemma 2.5.1]{MR2383890}).
Furthermore, by using the Weil descent datum $\alpha$, we can define the actions of $W_F$ on
$H^i_{\mathrm{LT}/\Gamma,c}$ and $H^i_{\mathrm{LT}/\Gamma}$ as follows.
For $w\in W_F$, let $\nu(w)$ denote the integer satisfying $w\vert_{F^\mathrm{ur}}=\sigma^{\nu(w)}$.
By taking the fiber product of diagrams
\[
 \xymatrix{%
 \Spa (\widehat{\overline{F}},\mathcal{O}_{\widehat{\overline{F}}})\ar[r]^-{w^*}\ar[d]&\Spa (\widehat{\overline{F}},\mathcal{O}_{\widehat{\overline{F}}})\ar[d]\\
 \Spa (\breve{F},\mathcal{O}_{\breve{F}})\ar[r]^-{(\sigma^*)^{\nu(w)}}&\Spa (\breve{F},\mathcal{O}_{\breve{F}})\lefteqn{,}
 }
 \qquad
 \xymatrix{%
 M_K/\Gamma\ar[r]^-{\alpha^{\nu(w)}}\ar[d]&M_K/\Gamma\ar[d]\\
 \Spa (\breve{F},\mathcal{O}_{\breve{F}})\ar[r]^-{(\sigma^*)^{\nu(w)}}&\Spa (\breve{F},\mathcal{O}_{\breve{F}})\lefteqn{,}
 }
\]
we obtain an isomorphism $\alpha_w\colon (M_K/\Gamma)\otimes_{\breve{F}}\widehat{\overline{F}}\to (M_K/\Gamma)\otimes_{\breve{F}}\widehat{\overline{F}}$ of adic spaces. The action of $w$ is defined to be $\alpha_w^*$.
By these constructions, we obtain two representations $H^i_{\mathrm{LT}/\Gamma,c}$ and $H^i_{\mathrm{LT}/\Gamma}$
of $\GL_n(F)\times D^\times\times W_F$.

Recall that any admissible representation $V$ of $\GL_n(F)/\Gamma$ is decomposed canonically into
$V=(\bigoplus_{\pi}V_\pi)\oplus V_{\text{non-cusp}}$, where
\begin{itemize}
 \item $\pi$ runs through irreducible supercuspidal representations of $\GL_n(F)$
       whose central characters are trivial on $\Gamma$, 
 \item $V_\pi$ is a direct sum of finitely many copies of $\pi$,
 \item and $V_{\text{non-cusp}}$ has no supercuspidal subquotient
\end{itemize} 
(see \cite[1.11, Variantes c)]{MR771671}).
We call $V_\pi$ the $\pi$-isotypic component of $V$. 
By definition we have $(V_{\pi})^\vee=(V^\vee)_{\pi^\vee}$
and $(V_{\text{non-cusp}})^\vee=(V^\vee)_{\text{non-cusp}}$.

We fix an isomorphism $\overline{\Q}_\ell\cong \C$ and identify them.
Here is a form of the non-abelian Lubin-Tate theory.

\begin{thm}[\cite{MR1876802}, \cite{MR1719811}, \cite{non-cusp}]\label{thm:NALT}
 For an irreducible supercuspidal representation
 $\pi$ of $\GL_n(F)$ whose central character is trivial on $\Gamma$, we have
 \[
 H^{n-1}_{\mathrm{LT}/\Gamma,c,\pi}(\tfrac{n-1}{2})=\pi\boxtimes \JL(\pi)^\vee\boxtimes\rec_F(\pi)^\vee
 \]
 as representations of $\GL_n(F)\times D^\times\times W_F$.
 Here $(\tfrac{n-1}{2})$ denotes the twist by the character $W_F\to \C^\times; w\mapsto q^{\frac{n-1}{2}\nu(w)}$,
 and $\JL(\pi)$ denotes the irreducible smooth representation of $D^\times$ attached to $\pi$
 under the local Jacquet-Langlands correspondence.
 Unless $i=n-1$, we have $H^i_{\mathrm{LT}/\Gamma,c,\pi}=0$.
\end{thm}

The following theorem was obtained in \cite{non-cusp}, in the course of the proof of 
the latter part of Theorem \ref{thm:NALT}.

\begin{thm}\label{thm:non-cusp}
 For every integer $i$, the kernel and cokernel of the natural map 
 $H^i_{\mathrm{LT}/\Gamma,c}\to H^i_{\mathrm{LT}/\Gamma}$ have no supercuspidal subquotient as representations
 of $\GL_n(F)$.
 In particular, for every irreducible supercuspidal representation $\pi$ of $\GL_n(F)$ whose central
 character is trivial on $\Gamma$,
 the induced map $H^i_{\mathrm{LT}/\Gamma,c,\pi}\to H^i_{\mathrm{LT/\Gamma,\pi}}$
 is an isomorphism.
\end{thm}

\begin{defn}\label{defn:trace-map}
 For a compact open subgroup $K$ of $\GL_n(\mathcal{O}_F)$, put $\Tr_K=(\GL_n(\mathcal{O}_F):K)^{-1}\Tr_{M_K}$,
 where $\Tr_{M_K}$ denotes the trace map
 $H_c^{2(n-1)}((M_K/\Gamma)\otimes_{\breve{F}}\widehat{\overline{F}},\overline{\Q}_\ell)(n-1)\to \overline{\Q}_\ell$.
 It is easy to see that $\Tr_K$ is compatible with the change of $K$.
 We write $\Tr$ for the homomorphism $H^{2(n-1)}_{\mathrm{LT}/\Gamma,c}(n-1)\to \overline{\Q}_\ell$
 induced from $\{\Tr_K\}_K$.
\end{defn}

\begin{prop}\label{prop:cup-product}
 Let $\pi$ be an irreducible supercuspidal representation of $\GL_n(F)$
 whose central character is trivial on $\Gamma$.
 Then, the cup product pairing
 \[
 \Tr(-\cup -)\colon H^{n-1}_{\mathrm{LT}/\Gamma,c}(\tfrac{n-1}{2})\times H^{n-1}_{\mathrm{LT}/\Gamma,c}(\tfrac{n-1}{2})\to \overline{\Q}_\ell
 \]
 induces a $D^\times\times W_F$-invariant pairing 
 $H^{n-1}_{\mathrm{LT}/\Gamma,c,\pi^\vee}(\tfrac{n-1}{2})\times H^{n-1}_{\mathrm{LT}/\Gamma,c,\pi}(\tfrac{n-1}{2})\to \overline{\Q}_\ell$
 satisfying the following condition:
 \begin{quote}
  for every compact open subgroup $K$ of $\GL_n(F)$, the restriction of it to
  $(H^{n-1}_{\mathrm{LT}/\Gamma,c,\pi^\vee})^K(\tfrac{n-1}{2})\times (H^{n-1}_{\mathrm{LT}/\Gamma,c,\pi})^K(\tfrac{n-1}{2})$ is a perfect pairing.
 \end{quote}
\end{prop}

\begin{prf}
 First, by the Poincar\'e duality for $M_K/\Gamma$, we know that
 the cup product pairing $(H^{n-1}_{\mathrm{LT}/\Gamma,c})^K(\tfrac{n-1}{2})\times (H^{n-1}_{\mathrm{LT}/\Gamma})^K(\tfrac{n-1}{2})\to \overline{\Q}_\ell$ is perfect for every compact open subgroup $K$ of $\GL_n(\mathcal{O}_F)$.
 This tells us that the induced map
 \[
  H^{n-1}_{\mathrm{LT}/\Gamma}(\tfrac{n-1}{2})\to (H^{n-1}_{\mathrm{LT}/\Gamma,c}(\tfrac{n-1}{2}))^\vee
 \]
 is an isomorphism.
 By taking $\pi$-isotypic parts and composing with the isomorphism in Theorem \ref{thm:non-cusp},
 we obtain an isomorphism
 \[
  H^{n-1}_{\mathrm{LT}/\Gamma,c,\pi}(\tfrac{n-1}{2})\xrightarrow{\cong} H^{n-1}_{\mathrm{LT}/\Gamma,\pi}(\tfrac{n-1}{2})\xrightarrow{\cong} (H^{n-1}_{\mathrm{LT}/\Gamma,c,\pi^\vee}(\tfrac{n-1}{2}))^\vee.
 \]
 Therefore, for every compact open subgroup $K$ of $\GL_n(F)$, we have an isomorphism
 \[
  (H^{n-1}_{\mathrm{LT}/\Gamma,c,\pi})^K(\tfrac{n-1}{2})\xrightarrow{\cong} ((H^{n-1}_{\mathrm{LT}/\Gamma,c,\pi^\vee})^K(\tfrac{n-1}{2}))^\vee.
 \]
 It is easy to see that this isomorphism is induced from the restriction of the cup product pairing
 to $(H^{n-1}_{\mathrm{LT}/\Gamma,c,\pi^\vee})^K(\tfrac{n-1}{2})\times (H^{n-1}_{\mathrm{LT}/\Gamma,c,\pi})^K(\tfrac{n-1}{2})$.
 This concludes the proof.
\end{prf}

\subsection{Twisting operator}\label{subsec:twisting}
Here we use the notation introduced in the beginning of Section \ref{subsec:div-alg}.
We will construct the twisting operator $\theta\colon \mathcal{M}_m\to \mathcal{M}_m$.

First we consider the case where $F/F^+$ is an unramified quadratic extension.
In this case we have $F^{\mathrm{ur}}=(F^+)^{\mathrm{ur}}$.
We write $\tau$ for the unique element of $\Gal(F^{\mathrm{ur}}/F^+)$
lifting the $q'$th power Frobenius automorphism $\overline{\tau}$ on $\overline{\F}_q=\overline{\F}_{q'}$.
It extends $\tau\in\Gal(F/F^+)$, and satisfies $\tau^2=\sigma$.
For an object $S$ of $\mathbf{Nilp}$, we write $S^\tau$ for the object
$S\to \Spec\mathcal{O}_{\breve{F}}\xrightarrow{\tau^*}\Spec\mathcal{O}_{\breve{F}}$ of $\mathbf{Nilp}$.

We write $\overline{\tau}_*\mathbb{X}$ for the pull-back of the formal $\mathcal{O}_F$-module $\mathbb{X}$
by $\overline{\tau}^*\colon \Spec\overline{\F}_q\to \Spec\overline{\F}_q$.
On the other hand, we denote by $\mathbb{X}^\tau$ the formal group $\mathbb{X}$ endowed with
the $\mathcal{O}_F$-action twisted by $\tau$
(that is, $\mathcal{O}_F\xrightarrow{\tau}\mathcal{O}_F\to \End(\mathbb{X})$).
It is easy to see that $\overline{\tau}_*\mathbb{X}$ and $\mathbb{X}^\tau$ are
one-dimensional formal $\mathcal{O}_F$-modules of $\mathcal{O}_F$-height $n$
over $(\Spec \overline{\F}_q)^\tau\in\mathbf{Nilp}$.
Hence these are isomorphic as formal $\mathcal{O}_F$-modules.
We fix an isomorphism $\iota\colon \overline{\tau}_*\mathbb{X}\xrightarrow{\cong}\mathbb{X}^\tau$ between them.
This isomorphism induces an automorphism on $D$:

\begin{defn}\label{defn:aut-D-unram}
 \begin{enumerate}
  \item An element $h\in D=\End_{\mathcal{O}_F}(\mathbb{X})\otimes_\Z\Q$ determines an element
	$\tau(h)=\iota\circ \overline{\tau}_*h\circ\iota^{-1}\in \End_{\mathcal{O}_F}(\mathbb{X}^\tau)\otimes_\Z\Q=\End_{\mathcal{O}_F}(\mathbb{X})\otimes_\Z\Q=D$. This gives an isomorphism $\tau\colon D\to D$ such that 
	$\tau\vert_F=\tau\in\Gal(F/F^+)$.
  \item  We denote the composite of $\mathcal{O}_F$-isogenies
	 \[
	 \mathbb{X}\xrightarrow{\Frob_{\mathbb{X}}}\overline{\sigma}_*\mathbb{X}=\overline{\tau}^2_*\mathbb{X}\xrightarrow{\overline{\tau}_*\iota}\overline{\tau}_*\mathbb{X}^\tau\xrightarrow{\iota}\mathbb{X}
	 \]
	 by $t$. It is an element of $D^\times$.
 \end{enumerate}
\end{defn}

\begin{lem}\label{lem:t-unram}
 The element $t\in D^\times$ satisfies $\tau^2=\Int(t)$ and $\tau(t)=t$.
\end{lem}

\begin{prf}
 For $h\in D$, we have
 \begin{align*}
 \tau^2(h)&=\iota\circ\overline{\tau}_*(\iota\circ \overline{\tau}_*h\circ \iota^{-1})\circ\iota^{-1}
 =(\iota\circ\overline{\tau}_*\iota)\circ \overline{\sigma}_*h\circ (\iota\circ\overline{\tau}_*\iota)^{-1}\\
  &=(t\circ \Frob_{\mathbb{X}}^{-1})\circ \overline{\sigma}_*h\circ (t\circ \Frob_{\mathbb{X}}^{-1})^{-1}\\
  &=t\circ (\Frob_{\mathbb{X}}^{-1}\circ \overline{\sigma}_*h\circ \Frob_{\mathbb{X}})\circ t^{-1}.
 \end{align*}
 By the functoriality of the relative Frobenius morphisms, the following diagram
 is commutative:
 \[
  \xymatrix{%
 \mathbb{X}\ar[rr]^-{\Frob_{\mathbb{X}}}\ar[d]^-{h}&& \overline{\sigma}_*\mathbb{X}\ar[d]^-{\overline{\sigma}_*h}\\
 \mathbb{X}\ar[rr]^-{\Frob_{\mathbb{X}}}&& \overline{\sigma}_*\mathbb{X}\lefteqn{.}
 }
 \]
 Hence we have $\tau^2(h)=t\circ h\circ t^{-1}$, as desired.

 Next consider $\tau(t)$. We have
 \begin{align*}
  \tau(t)&=\iota\circ\overline{\tau}_*(\iota\circ\overline{\tau}_*\iota\circ \Frob_{\mathbb{X}})\circ \iota^{-1}
  =\iota\circ\overline{\tau}_*\iota\circ\overline{\sigma}_*\iota\circ \overline{\tau}_*\Frob_{\mathbb{X}}\circ \iota^{-1}\\
  &=\iota\circ\overline{\tau}_*\iota\circ\overline{\sigma}_*\iota\circ \Frob_{\overline{\tau}_*\mathbb{X}}\circ \iota^{-1}\stackrel{(*)}{=}\iota\circ\overline{\tau}_*\iota\circ\Frob_{\mathbb{X}}=t.
 \end{align*}
 The equality $(*)$ follows from the functoriality of the relative Frobenius morphisms with respect to 
 $\iota\colon \overline{\tau}_*\mathbb{X}\to \mathbb{X}$. This completes the proof.
\end{prf}

Now we construct an isomorphism $\theta\colon \mathcal{M}_m\to \mathcal{M}_m$
that makes the following diagram commute:
\[
 \xymatrix{%
 \mathcal{M}_m\ar[r]^-{\theta}\ar[d]& \mathcal{M}_m\ar[d]\\
 \Spf\mathcal{O}_{\breve{F}}\ar[r]^-{\tau^*}& \Spf\mathcal{O}_{\breve{F}}\lefteqn{.}
 }
\]

\begin{defn}\label{defn:theta-unram}
 Let $S$ be an object of $\mathbf{Nilp}$. For $(X,\rho,\eta)\in \mathcal{M}_m(S)$, we put $\theta(X,\rho,\eta)=(X^\tau,\rho\circ\phi_{\overline{S}}^*\iota,\eta^\tau)\in\mathcal{M}_m(S^\tau)$, where
 \begin{itemize}
  \item $X^\tau$ is the formal group $X$ over $S$ endowed with the $\mathcal{O}_F$-action twisted by $\tau$,
  \item $\rho\circ\phi_{\overline{S}}^*\iota$ is the $\mathcal{O}_F$-quasi-isogeny
	$\phi_{\overline{S^\tau}}^*\mathbb{X}=\phi_{\overline{S}}^*(\overline{\tau}_*\mathbb{X})\xrightarrow{\phi_{\overline{S}}^*\iota}\phi_{\overline{S}}^*\mathbb{X}^\tau\xrightarrow{\rho}X^\tau\times_S\overline{S}$,
  \item and $\eta^\tau$ is the composite of $(\mathcal{O}_F/\mathfrak{p}_F^m)^n\xrightarrow{\tau}(\mathcal{O}_F/\mathfrak{p}_F^m)^n$ and $\eta$.
 \end{itemize}
 This gives a bijection $\theta\colon \mathcal{M}_m(S)\xrightarrow{\cong}\mathcal{M}_m(S^\tau)$,
 and an isomorphism $\theta\colon \mathcal{M}_m\xrightarrow{\cong} \mathcal{M}_m$ which covers
 $\tau^*\colon \Spf \mathcal{O}_{\breve{F}}\to \Spf \mathcal{O}_{\breve{F}}$.

 The isomorphism $\theta$ is compatible with the transition maps of the tower $\{\mathcal{M}_m\}$.
 Hence it induces automorphisms of the towers $\{\mathcal{M}_m\}$ and $\{M_m\}$.
\end{defn}

\begin{lem}\label{lem:property-theta-unram}
 \begin{enumerate}
  \item For $g\in\GL_n(F)$, we have $g\circ\theta=\theta\circ \tau(g)$, where $\tau\colon \GL_n(F)\to \GL_n(F)$
	is the isomorphism in Example \ref{exa:GL_n}.
  \item For $h\in D^\times$, we have $h\circ\theta=\theta\circ \tau(h)$, where $\tau\colon D^\times\to D^\times$
	is the isomorphism in Definition \ref{defn:aut-D-unram} (i).
  \item We have $\theta^2=\alpha\circ t$ and $\alpha\circ \theta=\theta\circ\alpha$.
 \end{enumerate}
\end{lem}

\begin{prf}
 The claim (i) is clear from the definition of $\theta$.

 As for (ii), take $(X,\rho,\eta)\in\mathcal{M}_m(S)$. Then we have
 \[
  (h\circ\theta)(X,\rho,\eta)=(X^\tau,\rho\circ \phi_{\overline{S}}^*\iota\circ\phi_{\overline{S^\tau}}^*h,\eta^\tau).
 \]
 Since $\phi_{\overline{S}}^*\iota\circ\phi_{\overline{S^\tau}}^*h=\phi_{\overline{S}}^*(\iota\circ \overline{\tau}_*h)=\phi_{\overline{S}}^*(\tau(h)\circ\iota)=\phi_{\overline{S}}^*(\tau(h))\circ\phi_{\overline{S}}^*\iota$, we have
 \[
  (h\circ\theta)(X,\rho,\eta)=(X^\tau,\rho\circ \phi_{\overline{S}}^*(\tau(h))\circ\phi_{\overline{S}}^*\iota,\eta^\tau)=(\theta\circ\tau(h))(X,\rho,\eta).
 \]
 Thus $h\circ\theta=\theta\circ\tau(h)$, as desired.

 We prove (iii). 
 For $(X,\rho,\eta)\in\mathcal{M}_m(S)$, $\theta^2(X,\rho,\eta)$ equals
 $(X,\rho\circ\phi_{\overline{S}}^*\iota\circ \phi_{\overline{S^\tau}}^*\iota,\eta)$.
 Since $\phi_{\overline{S}}^*\iota\circ \phi_{\overline{S^\tau}}^*\iota=\phi_{\overline{S}}^*(\iota\circ \overline{\tau}_*\iota)=\phi_{\overline{S}}^*(t\circ \Frob_{\mathbb{X}}^{-1})=\phi_{\overline{S}}^*(t)\circ \phi_{\overline{S}}^*\Frob_{\mathbb{X}}^{-1}$, we have $\theta^2(X,\rho,\eta)=\alpha(t(X,\rho,\eta))$.
 Hence $\theta^2=\alpha\circ t$.
 Finally, by (ii) and Lemma \ref{lem:t-unram} we conclude that
 \[
  \alpha\circ\theta=\theta^2\circ t^{-1}\circ\theta=\theta^3\circ \tau(t)^{-1}=\theta^3\circ t^{-1}=\theta\circ \alpha.
 \]
\end{prf}

 We fix $c\in W_{F^+}$ such that $c\vert_{F^{\mathrm{ur}}}=\tau$. 
 Assume that the subgroup $\Gamma\subset F^\times$ is stable under $\tau$.
 Then, $\theta\colon M_m/\Gamma\to M_m/\Gamma$ is induced.
 By taking the fiber product of diagrams
 \[
 \xymatrix{%
 \Spa (\widehat{\overline{F}},\mathcal{O}_{\widehat{\overline{F}}})\ar[r]^-{c^*}\ar[d]&\Spa (\widehat{\overline{F}},\mathcal{O}_{\widehat{\overline{F}}})\ar[d]\\
 \Spa (\breve{F},\mathcal{O}_{\breve{F}})\ar[r]^-{\tau^*}&\Spa (\breve{F},\mathcal{O}_{\breve{F}})\lefteqn{,}
 }
 \qquad
 \xymatrix{%
 M_m/\Gamma\ar[r]^-{\theta}\ar[d]&M_m/\Gamma\ar[d]\\
 \Spa (\breve{F},\mathcal{O}_{\breve{F}})\ar[r]^-{\tau^*}&\Spa (\breve{F},\mathcal{O}_{\breve{F}})\lefteqn{,}
 }
\]
we obtain an isomorphism $\theta_c\colon (M_m/\Gamma)\otimes_{\breve{F}}\widehat{\overline{F}}\to (M_m/\Gamma)\otimes_{\breve{F}}\widehat{\overline{F}}$ of adic spaces.
It induces an automorphism $\theta_c^*$ on the cohomology $H^i_{\mathrm{LT}/\Gamma,c}$,
for which we simply write $\theta$.

\begin{cor}\label{cor:coh-theta-unram}
 The following equalities of automorphisms on $H^i_{\mathrm{LT}/\Gamma,c}$ hold.
 \begin{enumerate}
  \item For $g\in\GL_n(F)$, we have $\theta\circ g=\tau(g)\circ\theta$,
	where $\tau\colon \GL_n(F)\to \GL_n(F)$	is the isomorphism in Example \ref{exa:GL_n}.
  \item For $h\in D^\times$, we have $\theta\circ h=\tau(h)\circ\theta$, 
	where $\tau\colon D^\times\to D^\times$	is the isomorphism in Definition \ref{defn:aut-D-unram} (i).
  \item We have $\theta^2=t\circ c^2$ and $\theta\circ w=cwc^{-1}\circ\theta$ for every $w\in W_F$.
 \end{enumerate}
\end{cor}

\begin{prf}
 The claims (i) and (ii) follow from Lemma \ref{lem:property-theta-unram} (i), (ii), respectively.
 For (iii), it suffices to show
 $\theta_c^2=\alpha_{c^2}\circ t$ and $\alpha_w\circ \theta_c=\theta_c\circ \alpha_{cwc^{-1}}$.
 These are consequences of Lemma \ref{lem:property-theta-unram} (iii), 
 the definitions of $\alpha_w$ and $\theta_c$, and the equality $\nu(cwc^{-1})=\nu(w)$.
\end{prf}

Next we consider the case where $F/F^+$ is a ramified quadratic extension
(here we do not need the tameness assumption).
We also write $\tau$ for the unique non-trivial element of $\Gal(F^{\mathrm{ur}}/(F^+)^{\mathrm{ur}})$.
It gives an extension of the original $\tau\in\Gal(F/F^+)$.
Note that $\tau^2=1$ and $\overline{\tau}=1$, where $\overline{\tau}$ denotes the automorphism of
the residue field $\overline{\F}_q$ of $\mathcal{O}_{F^{\mathrm{ur}}}$ induced by $\tau$.
Further, we have $\sigma\circ\tau=\tau\circ\sigma$ as automorphisms of $F^{\mathrm{ur}}$.
For an object $S$ of $\mathbf{Nilp}$, we write $S^\tau$ for the object
$S\to \Spec\mathcal{O}_{\breve{F}}\xrightarrow{\tau^*}\Spec\mathcal{O}_{\breve{F}}$ of $\mathbf{Nilp}$.

As in the unramified case, we fix an isomorphism $\iota\colon \mathbb{X}\xrightarrow{\cong}\mathbb{X}^\tau$
between formal $\mathcal{O}_F$-modules over $(\Spec\overline{\F}_q)^\tau=\Spec\overline{\F}_q\in\mathbf{Nilp}$.

\begin{defn}\label{defn:aut-D-ram}
 \begin{enumerate}
  \item An element $h\in D=\End_{\mathcal{O}_F}(\mathbb{X})\otimes_\Z\Q$ determines an element
	$\tau(h)=\iota\circ h\circ\iota^{-1}\in \End_{\mathcal{O}_F}(\mathbb{X}^\tau)\otimes_\Z\Q=\End_{\mathcal{O}_F}(\mathbb{X})\otimes_\Z\Q=D$. This gives an isomorphism $\tau\colon D\to D$ such that 
	$\tau\vert_F=\tau\in\Gal(F/F^+)$.
  \item We denote the composite $\mathbb{X}\xrightarrow{\iota}\mathbb{X}^\tau\xrightarrow{\iota}\mathbb{X}$
	by $t$. It is an element of $D^\times$.
 \end{enumerate}
\end{defn}

\begin{lem}\label{lem:t-ram}
 The element $t\in D^\times$ satisfies $\tau^2=\Int(t)$ and $\tau(t)=t$.
\end{lem}

\begin{prf}
 Clear from definition.
\end{prf}

Exactly in the same way, we can construct an isomorphism $\theta\colon \mathcal{M}_m\to \mathcal{M}_m$
that makes the following diagram commute:
\[
 \xymatrix{%
 \mathcal{M}_m\ar[r]^-{\theta}\ar[d]& \mathcal{M}_m\ar[d]\\
 \Spf\mathcal{O}_{\breve{F}}\ar[r]^-{\tau^*}& \Spf\mathcal{O}_{\breve{F}}\lefteqn{.}
 }
\]
It induces automorphisms of the towers $\{\mathcal{M}_m\}$ and $\{M_m\}$.

\begin{lem}\label{lem:property-theta-ram}
 \begin{enumerate}
  \item For $g\in\GL_n(F)$, we have $g\circ\theta=\theta\circ \tau(g)$, where $\tau\colon \GL_n(F)\to \GL_n(F)$
	is the isomorphism in Example \ref{exa:GL_n}.
  \item For $h\in D^\times$, we have $h\circ\theta=\theta\circ \tau(h)$, where $\tau\colon D^\times\to D^\times$
	is the isomorphism in Definition \ref{defn:aut-D-ram} (i).
  \item We have $\theta^2=t$ and $\alpha\circ \theta=\theta\circ\alpha$.
 \end{enumerate}
\end{lem}

\begin{prf}
 As in the proof of Lemma \ref{lem:property-theta-unram}, it suffices to show $\theta^2=t$.
 For an object $S$ of $\mathbf{Nilp}$ and $(X,\rho,\eta)\in\mathcal{M}_m(S)$, we have
 \[
  \theta^2(X,\rho,\eta)=(X,\rho\circ\phi^*_{\overline{S}}\iota\circ\phi^*_{\overline{S}}\iota,\eta)
 =(X,\rho\circ\phi^*_{\overline{S}}t,\eta)=t(X,\rho,\eta),
 \]
 as desired (note that $\overline{S^\tau}=\overline{S}$).
\end{prf}

We fix $c\in W_{F^+}$ such that $c\vert_{F^{\mathrm{ur}}}=\tau$.
Assume that $\Gamma\subset F^\times$ is stable under $\tau$.
As in the unramified case, we obtain an isomorphism
$\theta_c\colon (M_m/\Gamma)\otimes_{\breve{F}}\widehat{\overline{F}}\to (M_m/\Gamma)\otimes_{\breve{F}}\widehat{\overline{F}}$ of adic spaces.
It induces an automorphism $\theta_c^*$ on the cohomology $H^i_{\mathrm{LT}/\Gamma,c}$,
for which we simply write $\theta$.

\begin{cor}\label{cor:coh-theta-ram}
 The following equalities of automorphisms on $H^i_{\mathrm{LT}/\Gamma,c}$ hold.
 \begin{enumerate}
  \item For $g\in\GL_n(F)$, we have $\theta\circ g=\tau(g)\circ\theta$,
	where $\tau\colon \GL_n(F)\to \GL_n(F)$	is the isomorphism in Example \ref{exa:GL_n}.
  \item For $h\in D^\times$, we have $\theta\circ h=\tau(h)\circ\theta$, 
	where $\tau\colon D^\times\to D^\times$	is the isomorphism in Definition \ref{defn:aut-D-ram} (i).
  \item We have $\theta^2=t\circ c^2$ and $\theta\circ w=cwc^{-1}\circ\theta$ for every $w\in W_F$.
 \end{enumerate}
\end{cor}

\begin{prf}
 Similar as Corollary \ref{cor:coh-theta-unram}.
\end{prf}

Finally, consider the case $F=F^+$.

\begin{defn}\label{defn:aut-D-triv}
 We put $\tau=\id_{D^\times}$, $t=1\in D^\times$, $c=1\in W_{F^+}$ and $\theta=\id$ on $H^i_{\mathrm{LT}/\Gamma,c}$.
 Then the same statements as in Corollaries \ref{cor:coh-theta-unram}, \ref{cor:coh-theta-ram} 
 obviously hold.
\end{defn}

Now we return to a general separable extension $F/F^+$ with $[F:F^+]\le 2$.

\begin{lem}\label{lem:theta-cup-prod}
 Assume that $\Gamma\subset F^\times$ is stable under $\tau$.
 The cup product pairing
 \[
  \Tr(-\cup -)\colon H_{\mathrm{LT}/\Gamma,c}(\tfrac{n-1}{2})\times H_{\mathrm{LT}/\Gamma,c}(\tfrac{n-1}{2})\to \overline{\Q}_\ell
 \]
 in Proposition \ref{prop:cup-product} satisfies $\Tr(\theta x,\theta y)=q^{-\frac{n-1}{2}\nu(c^2)}\Tr(x\cup y)$.
\end{lem}

\begin{prf}
 Recall that the isomorphism 
 $\theta_c\colon (M_m/\Gamma)\otimes_{\breve{F}}\widehat{\overline{F}}\to (M_m/\Gamma)\otimes_{\breve{F}}\widehat{\overline{F}}$ covers $c^*\colon \Spa(\widehat{\overline{F}},\mathcal{O}_{\widehat{\overline{F}}})\to \Spa(\widehat{\overline{F}},\mathcal{O}_{\widehat{\overline{F}}})$. 

 If $F/F^+$ is an unramified quadratic extension, $c$ induces the $q'$th power map 
 on $\mu_{\ell^k}(\widehat{\overline{F}})=\mu_{\ell^k}(\breve{F}^+)$.
 Therefore we have $\Tr(\theta x,\theta y)=q'^{-(n-1)}\Tr(x\cup y)$.
 Since $q=q'^2$ and $\nu(c^2)=1$, this equals $q^{-\frac{n-1}{2}\nu(c^2)}\Tr(x\cup y)$.
 
 Otherwise $c$ acts trivially on $\mu_{\ell^k}(\widehat{\overline{F}})=\mu_{\ell^k}(\breve{F}^+)$,
 and $\nu(c^2)=0$. Hence we have $\Tr(\theta x,\theta y)=\Tr(x\cup y)=q^{-\frac{n-1}{2}\nu(c^2)}\Tr(x\cup y)$.
\end{prf}

\begin{thm}\label{thm:parity-geom}
 Here we consider $(\tau,t)$ as in Definitions \ref{defn:aut-D-unram}, \ref{defn:aut-D-ram}, \ref{defn:aut-D-triv}.
 Let $\pi$ be an irreducible supercuspidal representation of $\GL_n(F)$
 which is conjugate self-dual under the setting in Example \ref{exa:GL_n}.
 Then, $\JL(\pi)$ is conjugate self-dual with respect to $\tau$. Moreover, we have
 \[
  C_{\rec_F(\pi)}=(-1)^{n-1}C_{\JL(\pi)},
 \]
 where $C_{\JL(\pi)}$ denotes the parity of $\JL(\pi)$ with respect to $(\tau,t)$.
\end{thm}

\begin{prf}
 Since $\pi$ is conjugate self-dual, its central character $\omega_\pi$ satisfies
 $\omega_\pi(\tau(z))=\omega_\pi(z)^{-1}$ for every $z\in F^\times\subset \GL_n(F)$.
 Hence, for a uniformizer $\varpi'$ of $F^+$, we have $\omega_\pi(\varpi'^2)=1$.
 Put $\Gamma=\varpi'^{2\Z}\subset (F^+)^\times\subset F^\times$. It is a $\tau$-stable
 discrete cocompact subgroup of $F^\times$ on which $\omega_\pi$ is trivial.

 Let $\tau\colon \GL_n(F)\to \GL_n(F)$ be as in Example \ref{exa:GL_n},
 $\tau=\Int(c)\colon W_F\to W_F$ as in Example \ref{exa:Weil-group}, and
 $\tau=(\tau,\tau,\tau)\colon \GL_n(F)\times D^\times\times W_F\to \GL_n(F)\times D^\times\times W_F$.
 Then, Corollaries \ref{cor:coh-theta-unram}, \ref{cor:coh-theta-ram} tell us that
 $\theta$ gives an isomorphism $H^{n-1}_{\mathrm{LT}/\Gamma,c}\xrightarrow{\cong}(H^{n-1}_{\mathrm{LT}/\Gamma,c})^\tau$.
 Since the character $W_F\to \C^\times$; $w\mapsto q^{\frac{n-1}{2}\nu(w)}$ is $\tau$-invariant,
 we have $H^{n-1}_{\mathrm{LT}/\Gamma,c}(\frac{n-1}{2})\xrightarrow{\cong}(H^{n-1}_{\mathrm{LT}/\Gamma,c}(\frac{n-1}{2}))^\tau$ by twisting.
 By taking $\pi^\vee$-isotypic parts and using $\pi^\tau=\pi^\vee$, we obtain an isomorphism
 $\theta\colon H^{n-1}_{\mathrm{LT}/\Gamma,c,\pi^\vee}(\frac{n-1}{2})\xrightarrow{\cong}(H^{n-1}_{\mathrm{LT}/\Gamma,c,\pi}(\frac{n-1}{2}))^\tau$
 of representations of $\GL_n(F)\times D^\times\times W_F$.

 Take a $\tau$-stable compact open subgroup $K$ of $\GL_n(F)$. Then, $\theta$ induces an isomorphism
 $(H^{n-1}_{\mathrm{LT}/\Gamma,c,\pi^\vee}(\frac{n-1}{2}))^K\xrightarrow{\cong}((H^{n-1}_{\mathrm{LT}/\Gamma,c,\pi}(\frac{n-1}{2}))^K)^\tau$
 of representations of $D^\times\times W_F$. Consider the pairing
 \begin{align*}
  &\langle\ ,\ \rangle\colon (H^{n-1}_{\mathrm{LT}/\Gamma,c,\pi}(\tfrac{n-1}{2}))^K\times (H^{n-1}_{\mathrm{LT}/\Gamma,c,\pi}(\tfrac{n-1}{2}))^K\\
  &\qquad\qquad \xrightarrow[\cong]{\theta^{-1}\times \id}(H^{n-1}_{\mathrm{LT}/\Gamma,c,\pi^\vee}(\tfrac{n-1}{2}))^K\times (H^{n-1}_{\mathrm{LT}/\Gamma,c,\pi}(\tfrac{n-1}{2}))^K\xrightarrow{\Tr(-\cup -)}\overline{\Q}_\ell.
 \end{align*}
 It satisfies $\langle (\tau(h),\tau(w))x,(h,w)y\rangle=\langle x,y\rangle$ for every $h\in D^\times$ and $w\in W_F$.
 Moreover, Proposition \ref{prop:cup-product} tells us that it is a perfect pairing.
 We have
 \begin{align*}
  \langle y,x\rangle&=\Tr(\theta^{-1}(y)\cup x)=(-1)^{n-1}\Tr(x\cup \theta^{-1}(y))\stackrel{(1)}{=}(-1)^{n-1}q^{\frac{n-1}{2}\nu(c^2)}\Tr(\theta(x)\cup y)\\
 &=(-1)^{n-1}q^{\frac{n-1}{2}\nu(c^2)}\langle \theta^2(x),y\rangle
  \stackrel{(2)}{=}(-1)^{n-1}q^{\frac{n-1}{2}\nu(c^2)}\langle q^{-\frac{n-1}{2}\nu(c^2)}(t,c^2)(x),y\rangle\\
  &=(-1)^{n-1}\langle (t,c^2)x,y\rangle.
 \end{align*}
 Here (1) follows from Lemma \ref{lem:theta-cup-prod}, and (2) from the identity $\theta^2=t\circ c^2$ on 
 $H^{n-1}_{\mathrm{LT}/\Gamma,c,\pi}$
 (Corollary \ref{cor:coh-theta-unram} (iii) and Corollary \ref{cor:coh-theta-ram} (iii)); 
 the factor $q^{-\frac{n-1}{2}\nu(c^2)}$ arises from the twist $(\frac{n-1}{2})$.

 Now we specify $K$.
 Since $\pi$ is supercuspidal, it is generic. 
 Hence by \cite[\S 5, Th\'eor\`eme]{MR620708}, there exists an integer $m\ge 0$ such that
 $\dim \pi^{K_1(m)}=1$. Here $K_1(m)$ is the subgroup of $\GL_n(\mathcal{O}_F)$ consisting of
 matrices $(g_{ij})$ with $g_{n,1},\ldots,g_{n,n-1}\in \mathfrak{p}_F^m$ and 
 $g_{n,n}\in 1+\mathfrak{p}_F^m$. Clearly $K_1(m)$ is $\tau$-stable. We take $K$ as $K_1(m)$.
 Then, Theorem \ref{thm:NALT} tells us that 
 $(H^{n-1}_{\mathrm{LT}/\Gamma,c,\pi}(\frac{n-1}{2}))^K\cong \JL(\pi)^\vee\boxtimes \rec_F(\pi)^\vee$
 as representations of $D^\times\times W_F$. 
 Since $(\pi^\vee)^K=(\pi^K)^\vee$ is also one-dimensional,
 the existence of 
 $\theta\colon (H^{n-1}_{\mathrm{LT}/\Gamma,c,\pi^\vee}(\frac{n-1}{2}))^K\xrightarrow{\cong}((H^{n-1}_{\mathrm{LT}/\Gamma,c,\pi}(\frac{n-1}{2}))^K)^\tau$ tells us that
 \[
  \JL(\pi)\boxtimes \rec_F(\pi)=\JL(\pi^\vee)^\vee\boxtimes \rec_F(\pi^\vee)^\vee\cong \JL(\pi)^{\vee\tau}\boxtimes \rec_F(\pi)^{\vee\tau}.
 \]
 Thus $\JL(\pi)$ is conjugate self-dual with respect to $\tau$.
 Finally, by the existence of the pairing $\langle\ ,\ \rangle$, we conclude that
 the parity of the irreducible representation $\JL(\pi)^\vee\boxtimes \rec_F(\pi)^\vee$ of $D^\times\times W_F$
 with respect to $(\tau\times \tau,(t,c^2))$ is equal to $(-1)^{n-1}$.
 Replacing $\pi$ by $\pi^\vee$, we get the same result for $\JL(\pi)\boxtimes \rec_F(\pi)$.
 Therefore, by Lemma \ref{lem:product} we have
 $C_{\JL(\pi)}C_{\rec_F(\pi)}=(-1)^{n-1}$, and $C_{\rec_F(\pi)}=(-1)^{n-1}C_{\JL(\pi)}$.
 This completes the proof.
\end{prf}

\subsection{Formal $\mathcal{O}_F$-module over $\overline{\F}_q$ and division algebra}\label{subsec:formal-division}
Our remaining task for proving Theorem \ref{thm:main} is to describe $(\tau,t)$ in 
Definitions \ref{defn:aut-D-unram} and \ref{defn:aut-D-ram} explicitly, under the assumption that
$F/F^+$ is at worst tamely ramified and quadratic.

First we consider the easier case where $F$ has equal characteristic.
In this case, we have $F=\F_q((\varpi))$. We can take a one-dimensional formal $\mathcal{O}_F$-module $\mathbb{X}$
over $\overline{\F}_q$ as follows:
\[
 [a]_{\mathbb{X}}(X)=aX\ (a\in\F_q),\quad [\varpi]_{\mathbb{X}}(X)=X^{q^n}.
\]
Any element $a\in\F_{q^n}$ gives an endomorphism $X\mapsto aX$ of $\mathbb{X}$. On the other hand, we write $\Pi$
for the endomorphism $X\mapsto X^q$ of $\mathbb{X}$.
Note that $\Pi a=a^q\Pi$ for $a\in\F_{q^n}$ and $\Pi^n=\varpi$ in $\End_{\mathcal{O}_F}(\mathbb{X})$.
These elements are known to generate $\End_{\mathcal{O}_F}(\mathbb{X})$, and
we have $\End_{\mathcal{O}_F}(\mathbb{X})=\F_{q^n}[\Pi]=\mathcal{O}_{F_n}[\Pi]$,
which is a maximal order of the central division algebra over $F$ with invariant $1/n$.

Assume that $F/F^+$ is an unramified quadratic extension. We may assume that $F^+=\F_{q'}((\varpi))$.
Then, $\overline{\tau}_*\mathbb{X}$ and $\mathbb{X}^\tau$ are described explicitly as follows:
\begin{align*}
 [a]_{\overline{\tau}_*\mathbb{X}}(X)&=\overline{\tau}(aX)=a^{q'}X\ (a\in\F_q),&
 [\varpi]_{\overline{\tau}_*\mathbb{X}}(X)&=\overline{\tau}(X^{q^n})=X^{q^n},\\
 [a]_{\mathbb{X}^\tau}(X)&=[\overline{\tau}(a)]_{\mathbb{X}}(X)=a^{q'}X\ (a\in\F_q),&
 [\varpi]_{\mathbb{X}^\tau}(X)&=X^{q^n}.
\end{align*}
Hence we may take $\iota=\id_{\mathbb{X}}\colon \overline{\tau}_*\mathbb{X}\xrightarrow{\cong}\mathbb{X}^\tau$.
The following lemma is immediate.

\begin{prop}\label{prop:eq-char-unram}
 The pair $(\tau,t)$ constructed from $\iota=\id_{\mathbb{X}}$ as in Definition \ref{defn:aut-D-unram}
 coincides with that in Definition \ref{defn:tau-on-D} (i).
\end{prop}

Next assume that $p\neq 2$ and $F/F^+$ is a ramified quadratic extension.
We may assume that $F^+=\F_q((\varpi^2))$. Then $\mathbb{X}^\tau$ is described as follows:
\[
 [a]_{\mathbb{X}^\tau}(X)=aX\ (a\in\F_q),\quad [\varpi]_{\mathbb{X}^\tau}(X)=[-\varpi]_{\mathbb{X}}(X)=-X^{q^n}.
\]
Take $\beta\in \overline{\F}_q$ such that $\beta^{q^n-1}=-1$, and put $\alpha=\beta^{q-1}$.
Then, we may take an isomorphism $\iota\colon \mathbb{X}\xrightarrow{\cong}\mathbb{X}^\tau$; $X\mapsto \beta^{-1}X$.

\begin{prop}\label{prop:eq-char-ram}
 The pair $(\tau,t)$ constructed from $\iota\colon X\mapsto \beta^{-1}X$ as in Definition \ref{defn:aut-D-ram}
 coincides with that in Definition \ref{defn:tau-on-D} (ii).
\end{prop}

\begin{prf}
 For $a\in\F_{q^n}\subset \mathcal{O}_{F_n}[\Pi]$, 
 we have $\tau(a)\colon X\mapsto \beta^{-1}a\beta X=aX$; that is, $\tau(a)=a$.
 On the other hand, we have $\tau(\Pi)\colon X\mapsto \beta^{-1}(\beta X)^q=\alpha X^q$,
 and thus $\tau(\Pi)=\alpha\Pi$. Clearly we have $t=\beta^{-2}$.
 Hence the pair $(\tau,t)$ coincides with that in Definition \ref{defn:tau-on-D} (ii).
\end{prf}

Now we consider the case where $F$ is a $p$-adic field.
We regard formal $\mathcal{O}_F$-modules over $\overline{\F}_q$ as $\varpi$-divisible $\mathcal{O}_F$-modules.
We use the Dieudonn\'e theory for $\varpi$-divisible $\mathcal{O}_F$-modules over $\overline{\F}_q$
developed in \cite[Chapitre I, \S B.8]{MR2441311}.
Here we identify $\mathcal{O}_{\breve{F}}$ with
$W_{\mathcal{O}_F}(\overline{\F}_q)=\mathcal{O}_F\otimes_{W(\F_q)}W(\overline{\F}_q)$.
Let $\mathbb{D}=\mathcal{O}_{\breve{F}}^n$ be a free $\mathcal{O}_{\breve{F}}$-module of rank $n$.
We define a $\sigma$-linear map $F\colon \mathbb{D}\to \mathbb{D}$ and
a $\sigma^{-1}$-linear map $V\colon \mathbb{D}\to \mathbb{D}$ by
\[
 F(e_i)=\begin{cases}
	 \varpi e_{i+1}& i\neq n,\\ e_1& i=n,
       \end{cases}
 \quad
 V(e_i)=\begin{cases}
	 e_{i-1}& i\neq 1,\\ \varpi e_n& i=1,
	\end{cases}
\]
where $(e_1,\ldots,e_n)$ denotes the standard basis of $\mathbb{D}$.
Then, by \cite[Chapitre I, Proposition B.8.2]{MR2441311}, we can find a $\varpi$-divisible $\mathcal{O}_F$-module
$\mathbb{X}$ of $\mathcal{O}_F$-height $n$ over $\overline{\F}_q$
satisfying $\mathbb{D}_{\mathcal{O}_F}(\mathbb{X})\cong (\mathbb{D},F,V)$.
Since $V$ is topologically nilpotent and $\dim_{\overline{\F}_q}\mathbb{D}/V\mathbb{D}=1$,
$\mathbb{X}$ is a one-dimensional formal $\mathcal{O}_F$-module.

Let $D=F_n[\Pi]$ be the central division algebra over $F$ with invariant $1/n$ as in Section \ref{subsec:div-alg},
and $\mathcal{O}_D=\mathcal{O}_{F_n}[\Pi]$ its maximal order.
We will construct a homomorphism $\mathcal{O}_D\to \End_{\mathcal{O}_F}(\mathbb{X})$.
First, any $a\in\mathcal{O}_{F_n}$ defines an $\mathcal{O}_{\breve{F}}$-linear endomorphism on $\mathbb{D}$
by $e_i\mapsto \sigma^i(a)e_i$.
Since it commutes with $F$ and $V$, it gives an element of $\End_{\mathcal{O}_F}(\mathbb{X})$.
Let $\Pi$ be the $\mathcal{O}_{\breve{F}}$-linear endomorphism on $\mathbb{D}$ such that
\[
 \Pi(e_i)=\begin{cases}
	 e_{i-1}& i\neq 1,\\ \varpi e_n& i=1.
	\end{cases}
\]
It also commutes with $F$ and $V$, and gives an element of $\End_{\mathcal{O}_F}(\mathbb{X})$.
It is immediate to observe that $\Pi a=\sigma(a)\Pi$ for $a\in \mathcal{O}_{F_n}$ and $\Pi^n=\varpi$
as endomorphisms of $\mathbb{D}$. Therefore we obtain a homomorphism 
$\mathcal{O}_D=\mathcal{O}_{F_n}[\Pi]\to \End_{\mathcal{O}_F}(\mathbb{X})$, which is in fact an isomorphism.
In the following, we identify $\mathcal{O}_D$ and $\End_{\mathcal{O}_F}(\mathbb{X})$ by this isomorphism.

We assume that $F/F^+$ is an unramified quadratic extension, and take $\varpi$ in $F^+$.
Recall that in this case $\tau$ also denotes the unique element of $\Gal(F^{\mathrm{ur}}/F^+)$
lifting the $q'$th power Frobenius automorphism $\overline{\tau}$ on $\overline{\F}_q=\overline{\F}_{q'}$.
We describe $\overline{\tau}_*\mathbb{X}$ and $\mathbb{X}^\tau$ by means of the Dieudonn\'e module
as follows.

\begin{prop}\label{prop:Dieudonne-unram}
 Let $\tau_W\colon W(\overline{\F}_q)\to W(\overline{\F}_q)$ denote the homomorphism induced from
 $\overline{\mathbb{F}}_q\xrightarrow{\overline{\tau}}\overline{\mathbb{F}}_q$.
 We also write $\tau_W$ for the composite
 $W(\F_q)\to W(\overline{\F}_q)\xrightarrow{\tau_W} W(\overline{\F}_q)$.
 Note that the Dieudonn\'e module of a formal $\mathcal{O}_F$-module over
 $(\Spec\overline{\F}_q)^\tau\in\mathbf{Nilp}$ is a free
 $\mathcal{O}_F\otimes_{W(\F_q),\tau_W}W(\overline{\F}_q)$-module endowed with $F$ and $V$.
 We identify $\mathcal{O}_F\otimes_{W(\F_q),\tau_W}W(\overline{\F}_q)$ with $\mathcal{O}_{\breve{F}}$
 by the isomorphism $\mathcal{O}_F\otimes_{W(\F_q),\tau_W}W(\overline{\F}_q)\xrightarrow[\cong]{\tau\otimes \id}\mathcal{O}_F\otimes_{W(\F_q)}W(\overline{\F}_q)=\mathcal{O}_{\breve{F}}$.
 \begin{enumerate}
  \item For a formal $\mathcal{O}_F$-module $\mathbb{Y}$ over $\overline{\F}_q$, 
	we have $\mathbb{D}_{\mathcal{O}_F}(\overline{\tau}_*\mathbb{Y})=\tau_*\mathbb{D}_{\mathcal{O}_F}(\mathbb{Y})$
	and $\mathbb{D}_{\mathcal{O}_F}(\mathbb{Y}^\tau)=\mathbb{D}_{\mathcal{O}_F}(\mathbb{Y})$,
	where $\tau_*$ denotes the base change by
	$\tau\colon \mathcal{O}_{\breve{F}}\to \mathcal{O}_{\breve{F}}$.
  \item For $\mathbb{X}$ introduced above, we have $\mathbb{D}_{\mathcal{O}_F}(\overline{\tau}_*\mathbb{X})\cong \mathbb{D}_{\mathcal{O}_F}(\mathbb{X})$. 
 \end{enumerate}
\end{prop}

\begin{prf}
 We prove (i).
 By functoriality we have $\mathbb{D}_{\mathcal{O}_F}(\overline{\tau}_*\mathbb{Y})=(\id\otimes\tau_W)_*\mathbb{D}_{\mathcal{O}_F}(\mathbb{Y})$.
 Under the identification $\mathcal{O}_F\otimes_{W(\F_q),\tau_W}W(\overline{\F}_q)=\mathcal{O}_{\breve{F}}$,
 this equals 
 \[
  (\tau\otimes\id)_*(\id\otimes\tau_W)_*\mathbb{D}_{\mathcal{O}_F}(\mathbb{Y})=\tau_*\mathbb{D}_{\mathcal{O}_F}(\mathbb{Y}).
 \]
 On the other hand, we have $\mathbb{D}_{\mathcal{O}_F}(\mathbb{Y}^\tau)=(\tau^{-1}\otimes \id)_*\mathbb{D}_{\mathcal{O}_F}(\mathbb{Y})$.
 Under the identification, this clearly corresponds to the $\mathcal{O}_{\breve{F}}$-module
 $\mathbb{D}_{\mathcal{O}_F}(\mathbb{Y})$.

 The assertion (ii) is clear from the definition of $\mathbb{X}$ and the identification
 \[
  \tau_*\mathbb{D}\cong \mathbb{D}; (x_1,\ldots,x_n)\mapsto (\tau(x_1),\ldots,\tau(x_n)),
 \]
 as $\tau(\varpi)=\varpi$.
\end{prf}

\begin{prop}\label{prop:mixed-char-unram}
 Let $\iota\colon \overline{\tau}_*\mathbb{X}\xrightarrow{\cong}\mathbb{X}^\tau$ 
 be the isomorphism that induces the isomorphism in Proposition \ref{prop:Dieudonne-unram} (ii)
 on the Dieudonn\'e modules.
 The pair $(\tau,t)$ constructed from this $\iota$ as in Definition \ref{defn:aut-D-ram}
 coincides with that in Definition \ref{defn:tau-on-D} (i).
\end{prop}

\begin{prf}
 The claim on $\tau$ is clear from the definition.
 We will prove $\iota\circ\overline{\tau}_*\iota\circ \Frob_{\mathbb{X}}=\Pi$.
 Recall that $\Frob_{\mathbb{X}}\colon \mathbb{X}\to \overline{\sigma}_*\mathbb{X}$ induces
 $V\colon \mathbb{D}_{\mathcal{O}_F}(\mathbb{X})\to \sigma_*\mathbb{D}_{\mathcal{O}_F}(\mathbb{X})=\mathbb{D}_{\mathcal{O}_F}(\overline{\sigma}_*\mathbb{X})$.
 On the other hand,
 the composite $\sigma_*\mathbb{D}=\mathbb{D}_{\mathcal{O}_F}(\overline{\sigma}_*\mathbb{X})\xrightarrow{\mathbb{D}(\overline{\tau}_*\iota)}\mathbb{D}_{\mathcal{O}_F}(\overline{\tau}_*\mathbb{X}^\tau)\xrightarrow{\mathbb{D}(\iota)}\mathbb{D}_{\mathcal{O}_F}(\mathbb{X})=\mathbb{D}$ is
 equal to $(x_1,\ldots,x_n)\mapsto (\sigma(x_1),\ldots,\sigma(x_n))$.
 Since $\Pi(e_i)=V(e_i)$ for every $i$, we conclude that
 $\iota\circ\overline{\tau}_*\iota\circ \Frob_{\mathbb{X}}=\Pi$.
\end{prf}

Next we assume that $p\neq 2$ and $F/F^+$ is a ramified quadratic extension, and take $\varpi$
so that $\tau(\varpi)=-\varpi$.
Recall that in this case $\tau$ also denotes the unique non-trivial element of
$\Gal(F^{\mathrm{ur}}/(F^+)^{\mathrm{ur}})$.

\begin{prop}\label{prop:Dieudonne-ram}
 \begin{enumerate}
  \item For a formal $\mathcal{O}_F$-module $\mathbb{Y}$ over $\overline{\F}_q$, 
	we have $\mathbb{D}_{\mathcal{O}_F}(\mathbb{Y}^\tau)=\tau^{-1}_*\mathbb{D}_{\mathcal{O}_F}(\mathbb{Y})$,
	where $\tau^{-1}_*$ denotes the base change by
	$\tau^{-1}\colon \mathcal{O}_{\breve{F}}\to \mathcal{O}_{\breve{F}}$.
	For every $\mathcal{O}_F$-homomorphism $h\colon \mathbb{Y}\to \mathbb{Y}'$ between
	formal $\mathcal{O}_F$-modules over $\overline{\F}_q$,
	the homomorphism $\mathbb{D}_{\mathcal{O}_F}(\mathbb{Y}^\tau)\to \mathbb{D}_{\mathcal{O}_F}(\mathbb{Y}'^\tau)$ induced by $h\colon \mathbb{Y}^\tau\to \mathbb{Y}'^\tau$
	coincides with $\tau^{-1}_*\mathbb{D}(h)$.
  \item For $\mathbb{X}$ introduced above, we have $\mathbb{D}_{\mathcal{O}_F}(\mathbb{X}^\tau)=(\mathbb{D},F',V')$,
	where $F'$ and $V'$ are determined by
	\[
	 F'(e_i)=\begin{cases}
		 -\varpi e_{i+1}& i\neq n,\\ e_1& i=n,
		\end{cases}
	\quad
	V'(e_i)=\begin{cases}
		e_{i-1}& i\neq 1,\\ -\varpi e_n& i=1.
	       \end{cases}
	\]
  \item For an element $h\in \End_{\mathcal{O}_F}(\mathbb{X})$, 
	regard $\mathbb{D}(h)\colon \mathbb{D}\to \mathbb{D}$ as a matrix 
	$(h_{ij})\in M_n(\mathcal{O}_{\breve{F}})$.
	Then, the homomorphism $\mathbb{D}=\mathbb{D}_{\mathcal{O}_F}(\mathbb{X}^\tau)\to \mathbb{D}_{\mathcal{O}_F}(\mathbb{X}^\tau)=\mathbb{D}$
	induced by $h\colon \mathbb{X}^\tau\to \mathbb{X}^\tau$	is given by the matrix $(\tau^{-1}(h_{ij}))$.
 \end{enumerate}
\end{prop}

\begin{prf}
 The first assertion is clear from functoriality. The second is obvious
 from the definition of $\mathbb{X}$ and the identification
 \[
  \mathbb{D}\cong \tau^{-1}_*\mathbb{D}; (x_1,\ldots,x_n)\mapsto (\tau(x_1),\ldots,\tau(x_n)),
 \]
 as $\tau(\varpi)=-\varpi$.
 Let $h$ and $(h_{ij})$ be as in (iii).
 Under the identification $\mathbb{D}\cong \tau^{-1}_*\mathbb{D}$ above, $\tau_*^{-1}\mathbb{D}(h)$
 corresponds to $(\tau^{-1}(h_{ij}))$. The third assertion immediately follows from this. 
\end{prf}

\begin{prop}\label{prop:mixed-char-ram}
 As in Definition \ref{defn:tau-on-D} (ii), we take $\beta\in \mathcal{O}_{\breve{F}}$ such that $\beta^{q^n-1}=-1$
 and put $\alpha=\beta^{q-1}$.
 Let $\iota\colon \mathbb{X}\xrightarrow{\cong}\mathbb{X}^\tau$ be the isomorphism
 such that the induced homomorphism $\mathbb{D}=\mathbb{D}_{\mathcal{O}_F}(\mathbb{X})\to \mathbb{D}_{\mathcal{O}_F}(\mathbb{X}^\tau)=\mathbb{D}$ is given by
 $e_i\mapsto \sigma^i(\beta)^{-1}e_i$. 

 Then, the pair $(\tau,t)$ constructed from this $\iota$ as in Definition \ref{defn:aut-D-ram}
 coincides with that in Definition \ref{defn:tau-on-D} (ii).
\end{prop}

\begin{prf}
 For $a\in\mathcal{O}_{F_n}\subset \mathcal{O}_D$, the composite
 $\mathbb{D}=\mathbb{D}_{\mathcal{O}_F}(\mathbb{X}^\tau)\xrightarrow{\mathbb{D}(\iota\circ a\circ \iota^{-1})}\mathbb{D}_{\mathcal{O}_F}(\mathbb{X}^\tau)=\mathbb{D}$
 maps $e_i$ to $\sigma^i(a)e_i$. Hence, by Proposition \ref{prop:Dieudonne-ram} (iii),
 $\iota\circ a\circ\iota^{-1}\in \End_{\mathcal{O}_F}(\mathbb{X}^\tau)$ corresponds to 
 $\tau(a)\in \mathcal{O}_{F_n}\subset \End_{\mathcal{O}_F}(\mathbb{X})$ under the identification 
 $\End_{\mathcal{O}_F}(\mathbb{X})=\End_{\mathcal{O}_F}(\mathbb{X}^\tau)$.
 Similarly, the composite
 $\mathbb{D}=\mathbb{D}_{\mathcal{O}_F}(\mathbb{X}^\tau)\xrightarrow{\mathbb{D}(\iota\circ \Pi\circ \iota^{-1})}\mathbb{D}_{\mathcal{O}_F}(\mathbb{X}^\tau)=\mathbb{D}$ maps $e_i$ to
 \[
  \begin{cases}
   \frac{\sigma^i(\beta)}{\sigma^{i-1}(\beta)}e_{i-1}& i\neq 1,\\ \frac{\sigma(\beta)}{\sigma^n(\beta)}\varpi e_n& i=1.
  \end{cases}
 \]
 Since $\beta\in \mu_{2(q^n-1)}(\mathcal{O}_{\breve{F}})$, we have $\sigma(\beta)/\beta=\beta^{q-1}=\alpha$.
 Hence $\sigma^i(\beta)/\sigma^{i-1}(\beta)$ equals $\sigma^{i-1}(\alpha)$.
 Similarly, we have
 $\sigma(\beta)/\sigma^n(\beta)=\beta^q/\beta^{q^n}=\beta^q/(-\beta)=-\beta^{q-1}=-\alpha$.
 Noting that $\alpha\in \mu_{q^n-1}(\mathcal{O}_{\breve{F}})\subset \mathcal{O}_{F_n^+}$ is fixed by $\tau$,
 we can conclude that $\iota\circ\Pi\circ\iota^{-1}\in \End_{\mathcal{O}_F}(\mathbb{X}^\tau)$ corresponds to 
 $\alpha\Pi\in \End_{\mathcal{O}_F}(\mathbb{X})$ under the identification 
 $\End_{\mathcal{O}_F}(\mathbb{X})=\End_{\mathcal{O}_F}(\mathbb{X}^\tau)$.

 We can observe that the composite $\mathbb{X}\xrightarrow{\iota}\mathbb{X}^\tau\xrightarrow{\iota}\mathbb{X}$ equals $\beta^{-2}$ in the same way, by using the fact that $\beta\in \mu_{2(q^n-1)}(\mathcal{O}_{\breve{F}})\subset \mathcal{O}_{(F^+)^{\mathrm{ur}}}$ is fixed by $\tau$.
\end{prf}

By Theorem \ref{thm:parity-geom} and Propositions \ref{prop:eq-char-unram}, \ref{prop:eq-char-ram}, \ref{prop:mixed-char-unram}, \ref{prop:mixed-char-ram}, we complete the proof of Theorem \ref{thm:main}.

\section{The case of simple supercuspidal representations}\label{sec:simple-supercuspidal}
\subsection{Conjugate self-dual simple supercuspidal representations}
Here we apply our main theorem to simple supercuspidal representations.
Let the notation be as in Section \ref{subsec:div-alg}.
We briefly recall the notion of simple supercuspidal representations of $\GL_n(F)$ and $D^\times$.
See \cite{MR2730575}, \cite{MR3164986}, \cite{Knightly-Li} and \cite{Imai-Tsushima-JL} for detail.
Throughout this section, we fix a non-trivial additive character $\psi\colon \F_q\to \C^\times$
which factors through $\Tr_{\F_q/\F_p}\colon \F_q\to \F_p$.

First consider the case of $\GL_n(F)$.
Let us denote by $\Iw$ the standard Iwahori subgroup of $\GL_n(F)$, namely, the subgroup
of $\GL_n(\mathcal{O}_F)$ consisting of matrices whose image in $\GL_n(\F_q)$ is upper triangular.
We write $\Iw_+$ for the pro-$p$ unipotent radical of $\Iw$; it consists of matrices in $\Iw$
whose diagonal entries lie in $1+\mathfrak{p}_F$. 
Each element $\zeta\in \F_q^\times$ gives rise to a character
\[
 \psi_\zeta\colon \Iw_+\to \C^\times;\quad (a_{ij})\mapsto \psi(\overline{a_{12}}+\overline{a_{23}}+\ldots+\overline{a_{n-1,n}}+\zeta^{-1}\overline{\varpi^{-1}a_{n1}}).
\]
Here we denote the image of $a\in\mathcal{O}_F$ in $\overline{\F}_q$ by $\overline{a}$.

Let $\varphi_\zeta$ denote the matrix
\[
 \begin{pmatrix}
 0& 1& 0& \cdots& 0\\
 0& 0& 1& \cdots& 0\\
 \vdots& \vdots& \vdots& & \vdots\\
 0& 0& 0& \cdots& 1\\
 \widetilde{\zeta}\varpi& 0& 0& \cdots& 0
 \end{pmatrix},
\]
where $\widetilde{\zeta}$ denotes the unique element of $\mu_{q-1}(F)$ lifting $\zeta$.
It normalizes $\Iw_+$.
Put $H_\zeta=\mathcal{O}_F^\times\varphi_\zeta^\Z\Iw_+$.
It is an open compact-mod-center subgroup of $\GL_n(F)$ (note that it contains
the center $F^\times$, since $\varphi_\zeta^n=\widetilde{\zeta}\varpi$).
We write $(\F_q^\times)^\vee$ for the set of characters $\F_q^\times\to \C^\times$.
For a triple $(\zeta,\chi,c)\in \F_q^\times\times (\F_q^\times)^\vee\times \C^\times$,
define the character $\Lambda_{\zeta,\chi,c}\colon H_\zeta\to \C^\times$ by
\[
 \Lambda_{\zeta,\chi,c}(x)=\chi(\overline{x})\ (x\in \mathcal{O}_F^\times),\quad \Lambda_{\zeta,\chi,c}(\varphi_\zeta)=c,\quad \Lambda_{\zeta,\chi,c}\vert_{\Iw_+}=\psi_\zeta.
\]
We put $\pi_{\zeta,\chi,c}=\cInd_{H_\zeta}^{\GL_n(F)}\Lambda_{\zeta,\chi,c}$, which turns out
to be an irreducible supercuspidal representation of $\GL_n(F)$.
A representation obtained in this way is called a simple supercuspidal representation
of $\GL_n(F)$.
For another triple $(\zeta',\chi',c')\in \F_q^\times\times (\F_q^\times)^\vee\times \C^\times$, 
one can prove that $\pi_{\zeta,\chi,c}\cong \pi_{\zeta',\chi',c'}$ if and only if
$(\zeta,\chi,c)=(\zeta',\chi',c')$ (see \cite[Proposition 1.2]{Imai-Tsushima-JL}).
Thus simple supercuspidal representations of $\GL_n(F)$ are parameterized by the set 
$\F_q^\times\times (\F_q^\times)^\vee\times \C^\times$.

\begin{rem}
 Note that $\pi_{\zeta,\chi,c}$ implicitly depends on the choice of the uniformizer $\varpi$ of $F$.
 Later we take it as in Definition \ref{defn:tau-on-D}.
\end{rem}

The contragredient of $\pi_{\zeta,\chi,c}$ can be computed as follows:

\begin{prop}\label{prop:ssc-contragredient}
 For $(\zeta,\chi,c)\in \F_q^\times\times (\F_q^\times)^\vee\times \C^\times$, we have
 $\pi_{\zeta,\chi,c}^\vee\cong \pi_{(-1)^n\zeta,\chi^{-1},\chi(-1)c^{-1}}$.
\end{prop}

\begin{prf}
 For $a=\mathrm{diag}(1,-1,\ldots,(-1)^{n-1})$, we have $a\varphi_\zeta a^{-1}=-\varphi_{(-1)^n\zeta}$.
 As $a$ normalizes $\mathcal{O}_F^\times\Iw_+$, we obtain $aH_\zeta a^{-1}=H_{(-1)^n\zeta}$.
 Moreover, we can directly check that 
 $\Lambda_{\zeta,\chi,c}(h)^{-1}=\Lambda_{(-1)^n\zeta,\chi^{-1},\chi(-1)c^{-1}}(aha^{-1})$
 for $h\in H_\zeta$.
 Therefore $a$ intertwines $(H_\zeta,\Lambda_{\zeta,\chi,c}^{-1})$ and 
 $(H_{(-1)^n\zeta},\Lambda_{(-1)^n\zeta,\chi^{-1},\chi(-1)c^{-1}})$.
 By the same way as in the proof of Proposition \ref{prop:parity-cInd}, we conclude that
 \[
  \pi_{(-1)^n\zeta,\chi^{-1},\chi(-1)c^{-1}}=\cInd_{H_{(-1)^n\zeta}}^{\GL_n(F)}\Lambda_{(-1)^n\zeta,\chi^{-1},\chi(-1)c^{-1}}\cong \cInd_{H_{\zeta}}^{\GL_n(F)}\Lambda_{\zeta,\chi,c}^{-1}\cong \pi_{\zeta,\chi,c}^\vee.
 \]
\end{prf}

\begin{cor}\label{cor:ssc-csd}
 Let $(\zeta,\chi,c)$ be an element of $\F_q^\times\times (\F_q^\times)^\vee\times \C^\times$.
 \begin{enumerate}
  \item If $F/F^+$ is an unramified quadratic extension and $\varpi\in F^+$, then $\pi_{\zeta,\chi,c}$ is
	conjugate self-dual with respect to $\tau$ if and only if $\tau(\zeta)=(-1)^n\zeta$, $\chi^\tau=\chi^{-1}$
	and $c^2=\chi(-1)$, where $\tau$ denotes the $q'$th power Frobenius automorphism on $\F_q$.
  \item If $p\neq 2$, $F/F^+$ is a ramified quadratic extension and $\varpi$ satisfies $\tau(\varpi)=-\varpi$,
	then $\pi_{\zeta,\chi,c}$ is conjugate self-dual with respect to $\tau$ if and only if $n$ is odd,
	$\chi^2=1$ and $c^2=\chi(-1)$.
  \item If $F=F^+$, then $\pi_{\zeta,\chi,c}$ is conjugate self-dual with respect to $\tau=\id$
	(that is, self-dual) if and only if $n$ is even, $\chi^2=1$ and $c^2=\chi(-1)$.
 \end{enumerate}
\end{cor}

\begin{prf}
 In the proof of Proposition \ref{prop:parity-cInd}, we obtained an isomorphism 
 $(\cInd_H^G\chi)^\tau\cong \cInd_{H^\tau}^G\chi^\tau$. 
 We can use it to determine $(\pi_{\zeta,\chi,c})^\tau$ in each case as follows:
 \begin{enumerate}
  \item $(\pi_{\zeta,\chi,c})^\tau\cong \pi_{\tau^{-1}(\zeta),\chi^\tau,c}=\pi_{\tau(\zeta),\chi^\tau,c}$
	(note that $\psi(\tau(x))=\psi(x)$ for $x\in \F_q$).
  \item $(\pi_{\zeta,\chi,c})^\tau\cong \pi_{-\zeta,\chi,c}$.
  \item $(\pi_{\zeta,\chi,c})^\tau=\pi_{\zeta,\chi,c}$.
 \end{enumerate}
 Together with Proposition \ref{prop:ssc-contragredient}, we conclude the proof.
\end{prf}

Next we consider the case of $D^\times$. Let $(\zeta,\chi,c)$ be an element of $\F_q^\times\times (\F_q^\times)^\vee\times \C^\times$. Take $\xi\in \F_{q^n}^\times$ such that $\Nr_{\F_{q^n}/\F_q}(\xi)=\zeta$, and write $b$ for
the unique element of $\mu_{q^n-1}(\mathcal{O}_{F_n})$ lifting $\zeta$. 
Note that $(b\Pi)^n=\Nr_{F_n/F}(b)\Pi^n=\widetilde{\zeta}\varpi$.

Put $H^D_\xi=\mathcal{O}_F^\times(b\Pi)^\Z(1+\Pi\mathcal{O}_D)$. It is an open compact-mod-center subgroup of
$D^\times$. We define the character $\Lambda^D_{\xi,\chi,c}\colon H^D_\xi\to\C^\times$ by
\[
  \Lambda^D_{\xi,\chi,c}(x)=\chi(\overline{x})\ (x\in \mathcal{O}_F^\times),\ \Lambda^D_{\xi,\chi,c}(b\Pi)=c,\ \Lambda^D_{\xi,\chi,c}(1+b\Pi d)=\psi(\Tr_{\F_{q^n}/\F_q}(\overline{d}))\ (d\in\mathcal{O}_D).
\]
Here, $\overline{d}$ denotes the image of $d$ under
$\mathcal{O}_D\twoheadrightarrow\mathcal{O}_D/\Pi\mathcal{O}_D\stackrel{\cong}{\leftarrow} \mathcal{O}_{F_n}/\mathfrak{p}_{F_n}=\F_{q^n}$.
We put $\pi^D_{\zeta,\chi,c}=\cInd_{H^D_\xi}^{D^\times}\Lambda^D_{\xi,\chi,c}$, which turns out
to be an irreducible smooth representation of $D^\times$ whose isomorphism class depends only on $(\zeta,\chi,c)$.
A representation of $D^\times$, which is automatically supercuspidal,
obtained in this way is called a simple supercuspidal representation of $D^\times$.

The following theorem is proved in \cite[Theorem 3.5]{Imai-Tsushima-JL}.

\begin{thm}\label{thm:ssc-JL}
 For $(\zeta,\chi,c)\in \F_q^\times\times (\F_q^\times)^\vee\times\C^\times$, we have
 $\JL(\pi_{\zeta,\chi,c})=\pi^D_{\zeta,\chi,(-1)^{n-1}c}$.
\end{thm}

\subsection{Computation of parity}
Here we compute the parity of $\rec_F(\pi_{\zeta,\chi,c})$ for a conjugate self-dual simple supercuspidal
representation $\pi_{\zeta,\chi,c}$. We use Proposition \ref{prop:parity-cInd} to compute
the parity of $\pi^D_{\zeta,\chi,(-1)^{n-1}c}$.

\begin{prop}\label{prop:parity-JL}
 Let $(\zeta,\chi,c)$ be an element of $\F_q^\times\times(\F_q^\times)^\vee\times\C^\times$ such that
 $\pi_{\zeta,\chi,c}$ is conjugate self-dual with respect to $\tau$ under the setting in Example \ref{exa:GL_n}.
 \begin{enumerate}
  \item Suppose that $F/F^+$ is an unramified quadratic extension and $\varpi\in F^+$.
	Let $\varepsilon$ be an element of $\F_q^\times$ satisfying $\varepsilon^{q'-1}=-1$.
	Then the parity of $\pi^D_{\zeta,\chi,(-1)^{n-1}c}$ is equal to $(-1)^{n-1}\chi(\varepsilon)c$.
  \item Suppose that $p\neq 2$, $F/F^+$ is a ramified quadratic extension and $\tau(\varpi)=-\varpi$.
	Then the parity of $\pi^D_{\zeta,\chi,(-1)^{n-1}c}$ is equal to
	\[
	 \begin{cases}
	  1& \text{if $\chi$ is trivial},\\
	  -1& \text{if $\chi$ is non-trivial}.
	 \end{cases}
	\]
  \item Suppose that $F=F^+$. Then, the parity of $\pi^D_{\zeta,\chi,(-1)^{n-1}c}$ is equal to
	\[
	 \begin{cases}
	  1& \text{if $\chi$ is trivial},\\
	  -1& \text{if $\chi$ is non-trivial}.
	 \end{cases}
	\]
 \end{enumerate}
\end{prop}

\begin{prf}
 For simplicity, we put $c'=(-1)^{n-1}c$ and $\Psi=\Lambda^D_{\zeta,\chi,(-1)^{n-1}c}$.
 In each case we will find $a\in \mu_{q^n-1}(\mathcal{O}_{F_n})\subset D^\times$ which intertwines
 $(H^D_{\xi},\Psi^{-1})$ and $((H^D_{\xi})^\tau,\Psi^{\tau})$.

 Consider the case (i). Corollary \ref{cor:ssc-csd} tells us that $\zeta^{q'-1}=(-1)^n$, $\chi^\tau=\chi^{-1}$
 and $c^2=\chi(-1)$. Therefore we have
 $(\varepsilon\xi)^{(1+q+\cdots+q^{n-1})(q'-1)}=(-1)^{1+q+\cdots+q^{n-1}}\zeta^{q'-1}=(-1)^n\cdot (-1)^n=1$.
 Hence there exists $\eta\in\F_{q^n}^\times$ satisfying $\eta^{q'+1}=\varepsilon\xi$.
 Let $a_0$ be the unique element of $\mu_{q^n-1}(\mathcal{O}_{F_n})$ lifting $\eta$
 and put $a=\tau^{-1}(a_0)$.
 Since $\eta^{1-q}\xi^{q'}=(\varepsilon\xi)^{1-q'}\xi^{q'}=-\xi$, we have $a_0^{1-q}b^{q'}=-b$.
 Thus $a_0\tau(b\Pi)a_0^{-1}=a_0b^{q'}a_0^{-q}\Pi=-b\Pi$ and 
 $a(b\Pi)a^{-1}=-\tau^{-1}(b\Pi)$. In particular we have $aH^D_\xi a^{-1}=(H^D_\xi)^\tau$.

 Let us prove that $\Psi(h)^{-1}=\Psi^\tau(aha^{-1})$ for every $h\in H^D_\xi$. If $h\in \mathcal{O}_F^\times$,
 we have $\Psi(h)^{-1}=\chi(\overline{h})^{-1}=\chi^\tau(\overline{h})=\Psi^\tau(aha^{-1})$, as
 $\chi^\tau=\chi^{-1}$.
 If $h=b\Pi$, we have $\Psi(h)^{-1}=c'^{-1}$ and $\Psi^{\tau}(aha^{-1})=\Psi(-b\Pi)=\chi(-1)c'$.
 These are equal since $c'^2=c^2=\chi(-1)$.
 If $h=1+b\Pi d\in 1+\Pi\mathcal{O}_D$, we have $\Psi(h)^{-1}=\psi(\Tr_{\F_{q^n}/\F_q}(\overline{d}))^{-1}$
 and $\Psi^{\tau}(aha^{-1})=\Psi^{\tau}(1+a(b\Pi)a^{-1}\cdot ada^{-1})=\Psi(1-b\Pi a_0\tau(d)a_0^{-1})=\psi(\Tr_{\F_{q^n}/\F_q}(\overline{d}{}^{q'}))^{-1}$. Since $\psi$ factors through $\Tr_{\F_q/\F_{q'}}$, they are equal.

 Therefore $a$ intertwines $(H^D_{\xi},\Psi^{-1})$ and $((H^D_{\xi})^\tau,\Psi^{\tau})$.
 In this case, the element $z$ in Proposition \ref{prop:parity-cInd} becomes
 $a_0\Pi \tau^{-1}(a_0)=(a_0^{q'+1}b^{-1})\cdot (b\Pi)$.
 Note that the reduction of $a_0^{q'+1}b^{-1}\in\mu_{q^n-1}(\mathcal{O}_{F_n})$
 is equal to $\eta^{q'+1}\xi^{-1}=\varepsilon\in \F_q^\times$, and thus $a_0^{q'+1}b^{-1}$ lies in 
 $\mathcal{O}_F^\times$.
 Therefore, by Proposition \ref{prop:parity-cInd} the parity of $\pi^D_{\zeta,\chi,c'}$ is equal to
 \[
  \Psi(a_0\Pi \tau^{-1}(a_0))=\Psi(a_0^{q'+1}b^{-1})\Psi(b\Pi)=\chi(\varepsilon)c'=(-1)^{n-1}\chi(\varepsilon)c,
 \]
 as desired.

 Consider the case (ii). Corollary \ref{cor:ssc-csd} tells us that $n$ is odd, $\chi^2=1$ and $c^2=\chi(-1)$.
 Fix $\varepsilon\in\F_q^\times\setminus (\F_q^\times)^2$.
 As in Remark \ref{rem:tau-on-D} (ii), we can take $\tau\colon D\to D$ so that $\tau(\Pi)=-\Pi$,
 and $t$ as the unique element of $\mu_{q-1}(\mathcal{O}_F)$ lifting $\varepsilon$.
 Since $b\in\mu_{q^n-1}(\mathcal{O}_{F_n})\subset F_n^+$, we have $\tau(b\Pi)=-b\Pi$,
 and thus $H^D_\xi=(H^D_\xi)^\tau$.
 By the similar computation as in (i), we can observe that $\Psi(h)^{-1}=\Psi^\tau(h)$ for every $h\in H^D_\xi$.
 Therefore $1$ intertwines $(H^D_{\xi},\Psi^{-1})$ and $((H^D_{\xi})^\tau,\Psi^{\tau})$, and $z=t$.
 By Proposition \ref{prop:parity-cInd} the parity of $\pi^D_{\zeta,\chi,c'}$ is equal to
 $\Psi(t)=\chi(\varepsilon)\in \{\pm 1\}$.
 Since $\chi^2=1$, $\chi(\varepsilon)=1$ if and only if $\chi$ is trivial.

 Finally consider the case (iii).
 Corollary \ref{cor:ssc-csd} tells us that $n$ is even, $\chi^2=1$ and $c^2=\chi(-1)$.
 Take $\varepsilon\in \F_{q^2}^\times$ such that $\varepsilon^{q-1}=-1$, and let $a$ be the unique element
 of $\mu_{q^2-1}(\mathcal{O}_{F_2})$ lifting $\varepsilon$. 
 Since $n$ is even, $a$ belongs to $\mu_{q^n-1}(\mathcal{O}_{F_n})$.
 We have $a(b\Pi)a^{-1}=a^{1-q}b\Pi=-b\Pi$. Therefore $a$ normalizes $H^D_\xi$.
 By the similar computation as in (i), we can observe that $\Psi(h)^{-1}=\Psi(aha^{-1})$
 for every $h\in H^D_\xi$. 
 Therefore $a$ intertwines $(H^D_{\xi},\Psi^{-1})$ and $(H^D_{\xi},\Psi)$, and 
 $z=a^2\in\mu_{q^n-1}(\mathcal{O}_{F_n})$.
 Since $(\varepsilon^2)^{q-1}=1$, the reduction $\varepsilon^2$ of $z$ lies in $\F_q^\times$,
 and thus $z$ lies in $\mu_{q-1}(\mathcal{O}_F)$.
 Hence, by Proposition \ref{prop:parity-cInd} the parity of $\pi^D_{\zeta,\chi,c'}$ is equal to
 $\Psi(z)=\chi(\varepsilon^2)$. As $\varepsilon^2\in\F_q^\times\setminus (\F_q^\times)^2$ and $\chi^2=1$, 
 $\chi(\varepsilon^2)=1$ if and only if $\chi$ is trivial.
 This completes the proof.
\end{prf}

\begin{cor}\label{cor:parity-rec}
 Let $(\zeta,\chi,c)$ be as in Proposition \ref{prop:parity-JL}.
 \begin{enumerate}
  \item Suppose that $F/F^+$ is an unramified quadratic extension and $\varpi\in F^+$.
	Let $\varepsilon$ be an element of $\F_q^\times$ satisfying $\varepsilon^{q'-1}=-1$.
	Then the parity of $\rec_F(\pi_{\zeta,\chi,c})$ is equal to $\chi(\varepsilon)c$.
  \item Suppose that $p\neq 2$, $F/F^+$ is a ramified quadratic extension and $\tau(\varpi)=-\varpi$.
	Then the parity of $\rec_F(\pi_{\zeta,\chi,c})$ is equal to
	\[
	 \begin{cases}
	  1& \text{if $\chi$ is trivial},\\
	  -1& \text{if $\chi$ is non-trivial}.
	 \end{cases}
	\]
  \item Suppose that $F=F^+$. Then the parity of $\rec_F(\pi_{\zeta,\chi,c})$ is equal to
	\[
	 \begin{cases}
	  -1& \text{if $\chi$ is trivial},\\
	  1& \text{if $\chi$ is non-trivial}.
	 \end{cases}
	\]
 \end{enumerate}
\end{cor}

\begin{prf}
 Clear from Theorem \ref{thm:main}, Theorem \ref{thm:ssc-JL} and Proposition \ref{prop:parity-JL}.
 Recall that in the case (ii) (resp.\ (iii)), $n$ is odd (resp.\ even).
\end{prf}

\begin{rem}
 By Corollary \ref{cor:parity-rec} (iii), if a simple supercuspidal representation $\pi$ of $\GL_{2n}(F)$
 is self-dual and has trivial central character, $\rec_F(\pi)$ is symplectic and $\pi$ comes from
 $\mathrm{SO}(2n+1)$ by the endoscopic lifting. It is a starting point of a recent work of Oi \cite{Oi-ssc-SO}.

 On the other hand, if $F$ has characteristic $0$ and $p\neq 2$, Corollary \ref{cor:parity-rec} (i)
 has been obtained in \cite{Oi-ssc-U} by using the endoscopic character relation.
\end{rem}

\def\cftil#1{\ifmmode\setbox7\hbox{$\accent"5E#1$}\else
  \setbox7\hbox{\accent"5E#1}\penalty 10000\relax\fi\raise 1\ht7
  \hbox{\lower1.15ex\hbox to 1\wd7{\hss\accent"7E\hss}}\penalty 10000
  \hskip-1\wd7\penalty 10000\box7}
  \def\cftil#1{\ifmmode\setbox7\hbox{$\accent"5E#1$}\else
  \setbox7\hbox{\accent"5E#1}\penalty 10000\relax\fi\raise 1\ht7
  \hbox{\lower1.15ex\hbox to 1\wd7{\hss\accent"7E\hss}}\penalty 10000
  \hskip-1\wd7\penalty 10000\box7}
  \def\cftil#1{\ifmmode\setbox7\hbox{$\accent"5E#1$}\else
  \setbox7\hbox{\accent"5E#1}\penalty 10000\relax\fi\raise 1\ht7
  \hbox{\lower1.15ex\hbox to 1\wd7{\hss\accent"7E\hss}}\penalty 10000
  \hskip-1\wd7\penalty 10000\box7}
  \def\cftil#1{\ifmmode\setbox7\hbox{$\accent"5E#1$}\else
  \setbox7\hbox{\accent"5E#1}\penalty 10000\relax\fi\raise 1\ht7
  \hbox{\lower1.15ex\hbox to 1\wd7{\hss\accent"7E\hss}}\penalty 10000
  \hskip-1\wd7\penalty 10000\box7} \def\cprime{$'$} \def\cprime{$'$}
  \newcommand{\dummy}[1]{}
\providecommand{\bysame}{\leavevmode\hbox to3em{\hrulefill}\thinspace}
\providecommand{\MR}{\relax\ifhmode\unskip\space\fi MR }
\providecommand{\MRhref}[2]{%
  \href{http://www.ams.org/mathscinet-getitem?mr=#1}{#2}
}
\providecommand{\href}[2]{#2}

\end{document}